\documentclass[11pt,a4paper]{article}
\usepackage{a4wide}
\usepackage{amsfonts}
\usepackage{amsmath}
\usepackage{theorem}
\usepackage{color}
\usepackage{multirow}
\usepackage{framed}
\usepackage{booktabs}
\usepackage{graphicx}
\usepackage{hyperref}
\usepackage{mathtools}
\mathtoolsset{showonlyrefs}

\newtheorem{theorem}{Theorem}[section]
\newtheorem{lemma}[theorem]{Lemma}
\newtheorem{remark}[theorem]{Remark}
\newtheorem{algorithm}[theorem]{Algorithm}
\newtheorem{example}[theorem]{Example}
\newtheorem{corollary}[theorem]{Corollary}
\newtheorem{proposition}[theorem]{Proposition}

\def\tT{{\mbox{\tiny{T}}}}
\def\argmin{\mathop{\rm argmin}}

\DeclareMathOperator{\SO}{SO}
\newcommand{\HH}{\mathcal{H}}

\newcommand{\N}{\mathbb{N}}

\DeclareMathOperator{\dist}{dist}

\DeclareMathOperator{\grid}{{\mathcal G}}
\DeclareMathOperator{\neighbor}{\mathcal N}

\newcommand{\dx}{\,\mathrm{d}}
\newcommand{\TGV}{\operatorname{TGV}}
\newcommand{\TV}{\operatorname{TV}}

\newcommand{\SPD}{\mathcal P}
\newcommand{\Sym}[1]{\operatorname{Sym}}
\DeclareMathOperator{\interior}{int}

\DeclareMathOperator{\prox}{prox}
\DeclareMathOperator{\dom}{dom}

\DeclareMathOperator{\sgn}{sgn}

\DeclareMathOperator{\grad}{grad}
\newcommand{\geodesic}[1]{\ensuremath{\gamma_{\overset{\frown}{#1}}}}
\newcommand{\dotgeodesic}[1]{\ensuremath{\dot \gamma_{\overset{\frown}{#1}}}}

\newcommand{\MM}{\ensuremath{\mathcal{M}}}%Manifold 
\newcommand{\geoPV}[1]{\ensuremath{\gamma_{{#1}}}}
\newcommand{\RR}{\ensuremath{\mathbb{R}}}%Real values 
\newcommand{\dotgeoPV}[1]{\ensuremath{\dot\gamma_{{#1}}}}
\newcommand{\geo}[1]{\ensuremath{\gamma_{\overset{\frown}{#1}}}} 
\title{Recent Advances in Denoising of Manifold-Valued Images}
\date{December 2018}
\author{Bergmann, R.\thanks{Research Group Numerical Mathematics 
(Partial Differential Equations), Faculty of Mathematics, TU Chemnitz, Chemnitz, Germany. \emph{ronny.bergmann@mathematik.tu-chemnitz.de}}
\and
Laus, F.\thanks{Image Analysis Group, Department of Mathematics, TU Kaiserslautern, Kaiserslautern, Germany. 
\emph{$\{$friederike.laus,persch,steidl$\}$@mathematik.tu-kl.de}} \and Persch, J.\footnotemark[2] 
\and Steidl, G.\footnotemark[2]\ $^{,}$\footnote{Fraunhofer ITWM, Fraunhofer-Platz 1, D-67663 Kaiserslautern, Germany}
}
\begin{document}

\maketitle

%-----------------------------------------------------------------------------
\begin{abstract}
	Modern signal and image acquisition systems are able to capture  data that is no longer real-valued, but may take values on a manifold. 
	However, whenever measurements are taken, no matter whether manifold-valued or not, there occur tiny
inaccuracies, which result in noisy data.
In this chapter, we review recent advances in denoising of manifold-valued signals and images,
where we restrict our attention to variational models and appropriate minimization algorithms.
The algorithms are either classical as the subgradient algorithm or generalizations of 
the half-quadratic minimization method, the cyclic proximal point algorithm,
and the Douglas-Rachford algorithm to manifolds.
An important aspect when dealing with real-world data is the practical implementation. 
Here several groups provide software and toolboxes as the Manifold Optimization (Manopt) package and 
the manifold-valued image restoration toolbox (MVIRT).
%2available open source at www.mathematik.uni-kl.de/imagepro/members/bergmann/mvirt/
%3available open source at manopt.org
%4available open source at mtex-toolbox.github.io/
\end{abstract}
%-----------------------------------------------------------------------------
%-----------------------------------------------------------------------------
\section{INTRODUCTION} \label{introduction}
%-----------------------------------------------------------------------------
The mathematical notion of a manifold dates back to 1828, when Carl Friedrich Gauss established an important invariance property 
of surfaces while proving his Theorema Eregium. 
In his habilitation lecture in 1854, Bernhard Riemann intrinsically extended Gauss's theory 
making manifolds independent of their embedding in higher dimensional spaces.
This is now called a Riemannian manifold. 
Nowadays modern signal and image acquisition methods are able to capture information that is no longer restricted to Euclidean spaces 
but can be manifold-valued. Sophisticated models for human color perception involve non-Euclidean settings.
Moreover, it is often advantageous to model information from large data as points on a certain manifold.
Here are some typical examples.

\begin{figure}
	\centering
	\begin{tabular}{ccc}
	 \includegraphics[width = 0.3\textwidth]{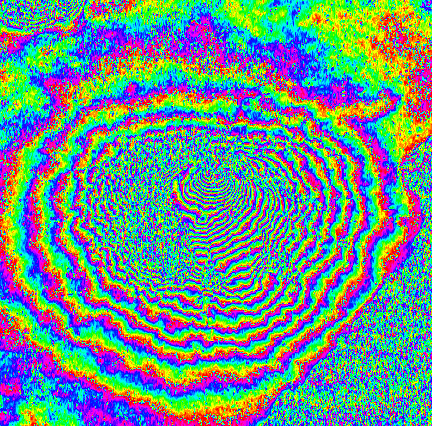}&
	 \includegraphics[width = 0.3\textwidth]{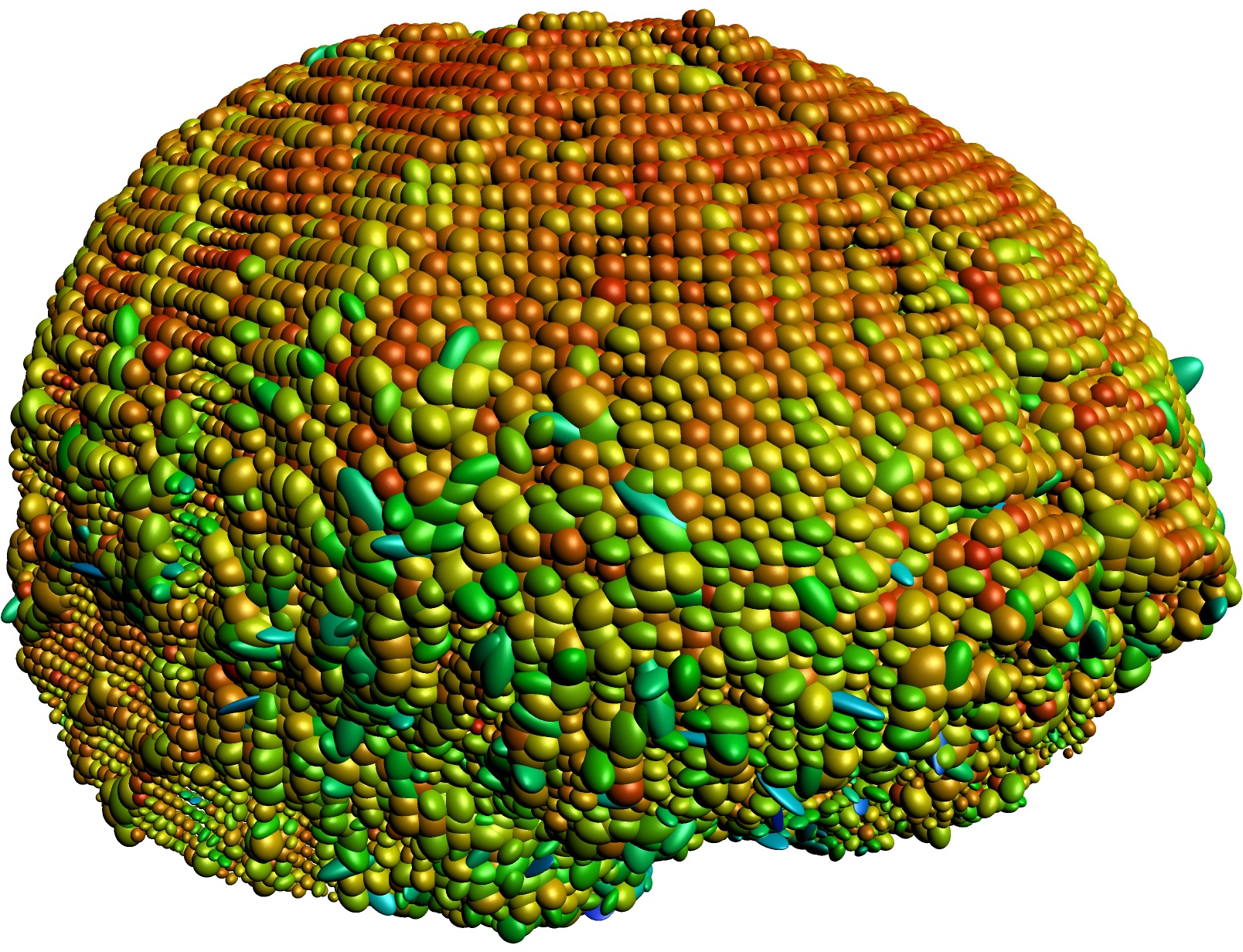}&
	 \includegraphics[width = 0.3\textwidth]{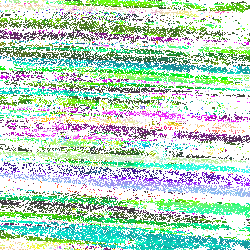}
	\end{tabular}
	\caption{Left: 
	InSAR data from Mt.~Vesuvius~\cite{RPG97}, 
	Middle: DT-MRI from the Camino project~\cite{camino} , 
	Right: EBSD data of an aluminum sample, Image courtesy: Institute of Materials Science and Engineering, TU Kaiserslautern }\label{fig:ex}
\end{figure}

\textbf{Interferometric Synthetic Aperture Radar (InSAR).}
	An InSAR image is obtained by calculating the phase difference of two Synthetic
	Aperture Radar (SAR) images of an area taken at different positions or times~\cite{BRF00}.
	It can be used to measure elevation when the measurements are taken at the same time or
	to detect millimeter-scale deformations over days or years. It has
	applications in the geophysical monitoring of natural hazards, for example earthquakes,
	volcanoes and landslides, and in structural engineering, in particular monitoring of subsidence and structural
	stability. InSAR produces phase-valued images, i.e., in each pixel the measurement lies on the
	circle $\mathbb S^1$, see Figure~\ref{fig:ex} (left).
	
\textbf{Image Color Spaces.} 
The usual RGB color space has a vector space structure, 
but there exist models which are physically better suited for the human color perception  as they encode luminance independent of color. 
Examples are the hue-saturation-value (HSV) color space or the related Ich space, where the hue
has its value on the circle $\mathbb S^1$, as well as the chromaticity-brightness (CB) color space,
where the chromaticity has values on the positive octant of the sphere $\mathbb S^2$.
For a recent geometric model of brightness perception we refer to \cite{BB2018}.
	
\textbf{Electron Backscatter Diffraction (EBSD).} 
EBSD is a microstructural crystallography characterization
	tech\-nique  used to analyze the microscopic
	structure of polycrystalline materials such as metals and minerals~\cite{BHS11}.
	Each point of a specimen is radiated by an electron beam and the diffraction
	pattern is measured, which gives information on the phase and the crystal
	orientation, a value on the rotation group $\SO(3)$. Regions of similar
	orientation are called grains. Material scientists are interested in the
	grain structure of the specimen, as it affects macroscopic properties such
	as ductility, electrical and lifetime properties. Since the atomic
	structure of a crystal is invariant under the symmetry of its atomic lattice, i.e., a symmetry group $S\subset\SO(3)$,
	the images obtained by EBSD have pixel values in $\SO(3)\slash  S$, see Figure \ref{fig:ex} (right). 
	The software MTEX \cite{MTEX} is designed for processing EBSD data.
	
\textbf{Diffusion Tensor Magnetic Resonance Imaging (DT-MRI).} 
In DT-MRI the diffusion of water molecules perpendicular to a magnetic field is measured in biological tissue. 
Taking at least six different data sets measured with different magnetic fields, 
a DT-MRI image with values in the manifold $\SPD(3)$ of symmetric positive definite $3\times 3$ matrices is computed at each pixel,
see Figure \ref{fig:ex} (middle).
The diffusion of water is influenced by the structure of the tissue.
Hence, the knowledge about the diffusion can be used to distinguish
between diseased and healthy parts. DT-MRI is a non-invasive and in-vivo technique,
which is extensively used in neurology, but can also be applied to detect defects
in other tissue, like muscles.

\textbf{Covariance Matrices in Texture Analysis and Brain Computer Interfaces.} 
Textures form a special class of images and appear at the same time as an important image feature.
There are different ways to model textures, one possibility is to encode the (local) dependence structure as a covariance matrix.
This has been used, e.g., in~\cite{TPM2006}, where textures are characterized based on the empirical covariance of certain features, 
for example color intensities and derivatives of different orders. 
Under the assumption that the empirical covariance matrices are non-degenerated we are again faced with the manifold $\SPD(d)$ of symmetric positive definite
$d \times d$ matrices,  where $d$ equals the number of features.
In Brain Computer Interfaces (BCI) approaches, EEG curves related to different human activities are measured at different brain positions at the same time.
Their covariance matrices can be used to analyze the corresponding activity, see, e.g. \cite{YCT2018}.

Beyond these applications, images with values in $\mathbb S^2$ appear when dealing
with 3D directional information~\cite{KS2002} as well as in the analysis of liquid crystals~\cite{Alo97}.
The rotation group $\operatorname{SO}(3)$ and the special Euclidean group $\operatorname{SE}(3)$ are considered in
tracking and (scene) motion analysis~\cite{rosman2012group,TPM08}.

Processing manifold-valued signals and images is a new challenge
that affects classical tasks like image restoration (denoising, inpainting), 
segmentation and clustering \cite{BFPS16}, 
registration and large deformation diffeomorphic mapping (LDDMM) or metamorphosis between different images, 
see, e.g., \cite{MuozMoreno2009ReviewOT,NPS17,ZNSY2014}. 

This chapter focuses on variational denoising methods for manifold-valued images.
The simplest idea is to embed the manifold into the Euclidean space and to apply the Euclidean models
with the constraint that the image values have to lie in a manifold.
Recall that by Whitney's theorem every smooth $d$-dimensional
manifold can be smoothly embedded into an Euclidean space of
dimension~$2d$. Such an approach was
given e.g., in \cite{RTKB14}. 
The advantage is that optimization algorithms in Euclidean spaces can be applied, 
where the models are in general non-convex due to the constraints.
So-called lifting schemes were proposed for circular-valued data 
in \cite{CS13} and for more general manifold-valued images in \cite{LSKC13}.
There, the problem is reformulated as a multilabel
optimization problem which is approached using convex relaxation techniques.
We also like to mention that current state-of-the-art methods for denoising of real-valued images 
include stochastic nonlocal patch-based approaches, 
as, e.g., the nonlocal Bayes' algorithm~\cite{LBM13}. 
A  generalization of this minimum mean square estimator (MMSE) based method 
to manifold-valued images was proposed in \cite{LNPS17}.

This chapter deals with spatially discrete, intrinsic  models and algorithms.
Note that in \cite{GMS93,GM07}, the notion of total variation of spatially continuous functions having their values on a
manifold was investigated where the authors apply the theory of Cartesian currents.
For a spatial continuous setting, the reader may also consult the recent paper \cite{CHS2019}.

An important aspect when working with real data is the practical implementation of the
developed methods and algorithms. 
In the spirit of reproducible research, several groups provide
their software and toolboxes, 
e.g., the Manifold Optimization (Manopt) package \cite{BMAS2014} and the manifold-valued image restoration toolbox (MVIRT~\cite{Bergmann2017}.

The outline of this chapter is as follows: Starting with the necessary preliminaries in Section~\ref{sec:prelim}, 
we review several denoising models for
manifold-valued images in Section~\ref{sec:models}. Appropriate minimization algorithms are discussed in Section \ref{sec:algs}.
Numerical examples in Section \ref{sec:numerics} illustrate the proposed methods.

%--------------------------------------------------------------------------
%-----------------------------------------------------------------------------
\section{Preliminaries on Riemannian Manifolds} \label{sec:prelim}
%-----------------------------------------------------------------------------
\subsection{General Notation}
%----------------------------------
Throughout this chapter, let $\mathcal M$ be a
connected, complete $d$-dimensional
Riemannian ma\-ni\-fold.
By \(T_x\mathcal M\) we denote the tangent space of $\mathcal M$
at \(x\in\mathcal M\) with the Riemannian
metric~$\langle\cdot,\cdot\rangle_x$ and corresponding norm~$\| \cdot \|_x$.
Further, let $T\MM$ be the tangent bundle of $\MM$. 
By
$\operatorname{dist} \colon \mathcal M \times \mathcal M
\rightarrow \mathbb R_{\ge 0}$
we denote the geodesic distance on \(\mathcal M\).
Let $\mathcal M^n$ be the product or \(n\)-fold power manifold
with product distance
\begin{equation}
\dist^2(x,y) \coloneqq \Big( \sum_{j=1}^n \operatorname{dist}^2(x_j,y_j) \Big)^\frac{1}{2}.
\end{equation}
Let~\(\gamma_{\overset{\frown}{x,y}}:[0,1] \rightarrow {\mathcal M}\) 
be a (not necessarily shortest)
geodesic connecting~\(x,y\in\mathcal M\).
Further, we apply the notation $\geoPV{x;\xi}$ to characterize the geodesics 
by its starting point $\geoPV{x;\xi}(0) = x$ and direction $\dotgeoPV{x;\xi}(0)=\xi \in T_x{\mathcal M}$.
Note that the geodesic $\gamma_{\overset{\frown}{x,y}}$ is unique
on manifolds with nonpositive
curvature. 
The exponential map $\exp_x\colon T_x\MM \to\MM$ is defined by 
\begin{equation}
\exp_x(\xi) \coloneqq \geoPV{x;\xi}(1).
\end{equation}
Since $\mathcal M$ is connected and complete, we know by the  Hopf-Rinow theorem that the exponential map
is indeed defined on the whole tangent space.
The exponential map realizes a local diffeomorphism 
from a neighborhood $\mathcal{D}_T(0_{x})$ of the origin $0_{x}$ of $T_{x}\MM$ 
into a neighborhood of $x \in \MM$.
More precisely,
extending the geodesic $\geoPV{x;\xi}$ from $t=0$ to infinity
is either minimizing $\dist(x,\geoPV{x;\xi}(t))$ all along or up to a finite time $t_0$ 
and not any longer afterwards.
In the latter case, $\geoPV{x;\xi}(t_0)$ is called cut point 
and the set of all cut points of all geodesics starting from $x$ 
is the cut locus $\mathcal{C}(x)$. 
This allows to define the inverse exponential map, also known as logarithmic map as
\begin{equation} \label{star}
\log_{x} \coloneqq \exp_{x}^{-1}\colon \MM \backslash {\mathcal C} (x) \to T_{x}\MM.
\end{equation}
Then, the Riemannian distance between $x,y\in\MM$, for $y\notin \mathcal{C}(x)$, can be written as 
\begin{equation} \label{log_dist}
\dist(x,y) = \langle \log_x(y),\log_x(y)\rangle_{x}^{\frac{1}{2}} = \lVert\log_x(y)\rVert_x.
\end{equation} 
Let $F\colon\MM\to\mathcal{N}$  be a smooth mapping between manifolds and $\xi \in T_x\MM$.
The  linear mapping 
\begin{equation}\label{eq:differential}
DF(x)\colon T_x\MM \to T_{F(x)}\mathcal{N},\quad \xi\mapsto DF(x)[\xi],
\end{equation}
is called differential of $F$ at $x \in \MM$. 
Let $F\colon\MM_1\to\MM_2$ and $G\colon\MM_2\to\MM_3$ 
be two smooth mappings. 
Then the differential of their concatenation $G\circ F$ applied to $\xi\in T_x\MM_1$ 
is given by the chain rule
\begin{equation}
D(G\circ F)(x)[\xi] = D G\bigl( F(x)\bigr)\bigl[DF(x)[\xi]\bigr].
\end{equation}
For a function $f\colon\MM\to\RR$, the Riemannian gradient $\grad_{\MM}$ is defined by
\begin{equation}\label{eq:graddiff}
\langle\grad_{\MM}f(x), \xi\rangle_x \coloneqq Df(x)[\xi]
\end{equation}
for all $\xi\in T_x\MM$.
A mapping $R_p\colon\MM \rightarrow \MM$ 
is called geodesic reflection at $x\in \MM$, if
\begin{equation} \label{def:refl}
R_x(x) = x \quad \text{and} \quad D\bigl(R_x(x)\bigr) = -I.
\end{equation}
For $\MM = \mathbb R^n$ we simply have $R_p (x) = 2p - x$.
A connected Riemannian manifold $\MM$ is (globally)
symmetric if the geodesic reflection at any point $x \in \MM$
is an isometry of $\MM$, i.e. 
$
	\dist\left(R_p(x),R_p(y) \right) = \dist(x,y)
$
for all $x,y \in \MM$.
All manifolds considered in this chapter are symmetric ones.

\begin{figure*}[t]
	\centering
	\includegraphics{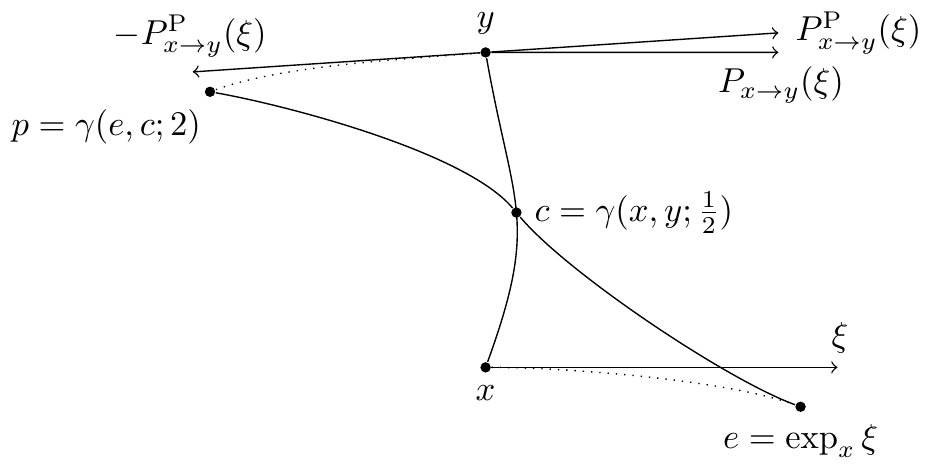}
	\hspace{0.2cm}
	\includegraphics{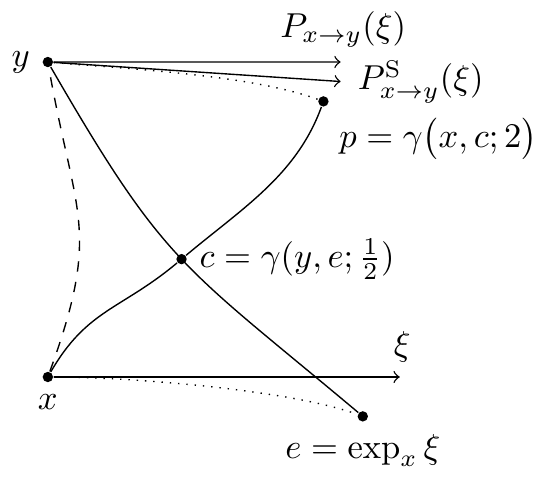}
	\caption{Illustration of pole ladder (left) and Schild's ladder (right) for the approximation of $P_{x\to y} \xi$.}\label{fig:pole+schild}
\end{figure*}

Let ${\mathcal X}(\MM)$ be the set of smooth vector fields on $\MM$.
Given a curve $\gamma\colon[0,1]\to\MM$, we denote by ${\mathcal X}(\gamma)$ the set of smooth vector fields along $\gamma$, 
i.e., $X\in{\mathcal X}(\gamma)$ is a smooth mapping $X\colon[0,1]\to T\MM$ with $X(t)\in T_{\gamma(t)}\MM$.
A vector field $X\in{\mathcal X}(\gamma)$ 
is called parallel to $\gamma\colon[0,1]\to\MM$,
if the covariant derivative along $\gamma$ fulfills
$\frac{D}{\dx t}X = 0$  for all $t\in[0,1]$.
We define the parallel transport of a tangent vector $\xi\in T_x\MM$ to $T_y\MM$  by 
\begin{equation}
P_{x\to y} \xi \coloneqq X(1),
\end{equation}
where $X\in\mathcal{X}(\geo{x,y})$ is the vector field parallel to a minimizing geodesic $\geo{x,y}$ with $X(0)=\xi$.
There exist analytical expressions of the parallel transport for few manifolds as spheres
or positive definite matrices. However, the parallel transport can be locally approximated
by Schild's ladder~\cite{EPS72,KMN2000}
or by the pole ladder~\cite{LP14Sch}.
Recently, it was shown that for connected, complete, symmetric Riemannian manifold, 
the pole ladder coincides with the parallel transport along geodesics~\cite{Pen18}.
Therefore, we prefer the pole ladder approach.
Given $x,y\in\MM$, 
the pole ladder transports $\xi\in T_x\MM$ to $T_y\MM$ by
\begin{equation}\label{pre:eq:pole}
P_{x\to y}^{\mathrm{P}}(\xi)\coloneqq -\log_{y}\Bigl(\gamma\Bigl(\exp_{x}(\xi),\gamma\bigl(x,y;\tfrac{1}{2}\bigr);2\Bigr)\Bigr)\in T_y\MM,\
\end{equation}
where we use the notation $\gamma(x,y;t) \coloneqq \gamma_{\overset{\frown}{x,y}} (t)$. 
For comparison, Schild's ladder transports as follows:
\begin{equation}\label{pre:eq:schild}
P_{x\to y}^{\mathrm{S}}(\xi)\coloneqq \log_{y}\Bigl(\gamma\Bigl(x,\gamma\bigl(y,\exp_{x}(\xi);\tfrac{1}{2}\bigr);2\Bigr)\Bigr)\in T_y\MM.
\end{equation}
Both transport schemes are illustrated in Figure~\ref{fig:pole+schild}.
%------------------------------------------------
%-----------------------------------------------
In our minimization algorithms, we will need the Riemannian gradient of special functions,
in particular of those appearing in the pole ladder~\eqref{pre:eq:pole}.
 These gradients can be derived from differentials, see~\eqref{eq:graddiff}, and can be computed for symmetric Riemannian manifolds using the theory of Jacobi fields.
The following lemma collects the final results which can be partially found 
in~\cite{BBSW16,BHSW17,Carmo92}. 
For the complete proof we refer to~\cite{Persch2018}.

\begin{lemma}\label{pre:differentials}
Let $\MM$ be a symmetric Riemannian ma\-ni\-fold and 
$F$ one of the functions i) - v) below with parameter $y$, resp.~$u$, together with 
the coefficient map $\alpha\colon \mathbb R \rightarrow \mathbb R$ and parameter $T$.
Then the differential $DF(x)$ at $x \in \MM$ is given for all $\xi \in T_x \MM$ by 
	\begin{equation}\label{alles}
	D F (x) [\xi] = \sum_{k=1}^d \langle \xi, \Xi_k(0)\rangle_{x} \alpha(\kappa_k) \Xi_k(T),
	\end{equation}
where $\{\Xi_k\}_{k=1}^d$ denotes a parallel transported orthogonal frame along the geodesic $\gamma$
with $\gamma(0) = x$ and $\gamma (1) = y$, where $y \coloneqq \exp_x(u)$ if $F$ depends on $u$.
Further, the frame diagonalizes the Riemannian curvature tensor 
$R(\cdot,\dot\gamma)\dot{\gamma}$ with respective eigenvalues $\kappa_k$, $k = 1,\dots,d$.
The functions $F$ and $\alpha$ are given as follows:
	\begin{enumerate}
		\item[i)] For $F\coloneqq \exp_{\cdot} (u)$, we have $T=1$,  $y \coloneqq \exp_{x} (u)$ and 
		\begin{equation}
		\alpha(\kappa) \coloneqq \begin{cases}
		\cosh(\sqrt{-\kappa})\quad&\kappa <0,\\		
		1\quad&\kappa =0,\\
		\cos(\sqrt{\kappa})\quad&\kappa >0.
		\end{cases}
		\end{equation}
		\item[ii)] \ For $F  \coloneqq \log_\cdot (y)$, we have $T=0$ and
		\begin{equation}
		\alpha(\kappa) \coloneqq \begin{cases}
		-\sqrt{-\kappa}\frac{\cosh(\sqrt{-\kappa})}{\sinh(\sqrt{-\kappa})}\quad&\kappa <0,\\		
		-1\quad&\kappa =0,\\
		-\sqrt{\kappa}\frac{\cos(\sqrt{\kappa})}{\sin(\sqrt{\kappa})}\quad&\kappa >0.
		\end{cases}
		\end{equation}
		
		\item[iii)]\ For $F \coloneqq \log_y (\cdot)$, we have  $T=1$ and
		\begin{equation}
		\alpha(\kappa) \coloneqq \begin{cases}
		\frac{\sqrt{-\kappa}}{\sinh(\sqrt{-\kappa})}\quad&\kappa <0,\\		
		1\quad&\kappa =0,\\
		\frac{\sqrt{\kappa}}{\sin(\sqrt{\kappa})}\quad&\kappa >0.
		\end{cases}
		\end{equation}
		\item[iv)] For $F\coloneqq \geodesic{\cdot,y}(\tau)$, we have $T=\tau$ and
		\begin{equation}
		\alpha(\kappa) \coloneqq \begin{cases}
		\frac{\sinh\bigl(\sqrt{-\kappa}(1-\tau)\bigr)}{\sinh(\sqrt{-\kappa})}\quad&\kappa <0,\\		
		1-\tau\quad&\kappa =0,\\
		\frac{\sin\bigl(\sqrt{\kappa}(1-\tau)\bigr)}{\sin(\sqrt{\kappa})}\quad&\kappa >0.
		\end{cases}
		\end{equation}
		\item[v)] For $F \coloneqq \geodesic{y,\cdot}(\tau)$, we have $T=1-\tau$ and
		\begin{equation}
		\alpha(\kappa) \coloneqq \begin{cases}
		\frac{\sinh(\sqrt{-\kappa}\tau)}{\sinh(\sqrt{-\kappa})}\quad&\kappa <0,\\		
		\tau\quad&\kappa =0,\\
		\frac{\sin(\sqrt{\kappa}\tau)}{\sin(\sqrt{\kappa})}\quad&\kappa >0.
		\end{cases}
		\end{equation}
\item[vi)]\ Finally, we obtain for $F \coloneqq 	\exp_x (\cdot)$  with
\begin{equation}
		\alpha(\kappa) = \begin{cases}
		\frac{\sinh(\sqrt{-\kappa})}{\sqrt{-\kappa}}\quad&\kappa <0,\\		
		1&\kappa =0,\\
		\frac{\sin(\sqrt{\kappa})}{\sqrt{\kappa}}\quad&\kappa >0,\\
		\end{cases}
		\end{equation}
		and $T=1$
		that the differential $D F (u)$ of $F$ at $u \in T_x \MM$ is given by \eqref{alles},
		where we have to replace $x \in \MM$ by $u \in T_x \MM$ and to set $y \coloneqq \exp_x(u)$.
\end{enumerate}
\end{lemma}

The adjoint operator $(D F)^* (x)\colon T_{F(x)}\MM\to T_x\MM$ of \eqref{alles}, which is also needed for the computation of the gradients, is given by
\begin{equation}\label{pre:eq:adjoint}
(D F)^* (x) [w] = \sum_{k=1}^d \langle w,\Xi_k\rangle_{F(x)} \alpha_k\xi_k, \quad w\in T_{F(x)}\MM.
\end{equation}

%---------------------------------------------------
\subsection{Convexity and Hadamard Manifolds}
%---------------------------------------------------
A subset $C\subseteq\MM$ is called weakly (strongly) convex if for all points $x,y\in C$ there exists a (unique) geodesic $\geodesic{x,y}$ of minimal length 
which is contained entirely in $C$.
Let $\varphi\colon\MM\supset C\to\RR\cup\{\infty\}$ be a real-valued function on a weakly convex set $C$, then $\varphi$ is called convex if 
\begin{equation}
\label{hada:conv}
\varphi\bigl(\geodesic{x,y}(t)\bigr) \le (1-t)\varphi(x)+t\varphi(y),
\end{equation}
for all $x,y\in C$. The function $\varphi$ is called strictly convex, if the above equation holds strictly for all $t\in(0,1)$. 
A function $\varphi\colon\MM\supseteq C\to\RR\cup\{\infty\}$ is $\kappa$-strongly convex if 
\begin{equation}\label{hada:strconv}
\varphi\bigl(\geodesic{x,y}(t)\bigr)\le (1-t)\varphi(x)+t\varphi(y)-\kappa t(1-t)\dist^2(x,y).
\end{equation}
Simply connected, complete Riemannian manifolds of nonpositive sectional curvature are called Hadamard manifolds.
We denote them by $\HH$.
Examples are the manifold of positive definite matrices of fixed size with the affine invariant metric or hyperbolic spaces.
An important property of Hadamard manifolds is, that the distance function $\dist(\cdot,\cdot)$ is jointly convex, which makes $\dist^2(\cdot,\cdot)$ strictly convex. Further, $\dist^2(\cdot,y)$ is a 1-strongly convex function on a Riemannian manifold if and only if $\HH$ is a Hadamard manifold.
The domain of a function $\varphi\colon\HH \to \RR\cup\{\infty\}$ is defined by
\begin{equation}
\dom(\varphi) \coloneqq\{x\in \HH:\varphi(x)< \infty\},
\end{equation}
in general we work with proper functions, i.e., 
$\dom(\varphi)\neq\emptyset$. 
A function $\varphi$ is called lower semi continuous (lsc) 
if the set $\{x \in \HH:\varphi(x)\le c\}$ is closed for all $c\in\RR$. 
Whenever $\dist(x,x_0)\to\infty$ for some $x_0\in\HH$, the function $\varphi$ is called coercive if $\varphi(x)\to\infty$. 
Concerning minimizers of convex functions, the next theorem summarizes some basic facts.
A proper, convex, lsc functions $\varphi\colon \HH \rightarrow \RR\cup\{\infty\}$ has a minimizer
	if it is coercive and unique minimizer if it is in addition strongly convex.	
For further information on more general Hadamard spaces and the basics of convex analysis therein, we refer to~\cite{Bac14}.

%----------------------------------------------------------------------------
%-----------------------------------------------------------------------------
\section{Intrinsic Variational Restoration Models} \label{sec:models}
%-----------------------------------------------------------------------------
%
We consider images as mappings from the image grid $\grid=\{1,\dots,n_1\}\times\{1,\dots,n_2\}$
to a  Riemannian manifold  $\MM$. Let $n \coloneqq n_1 n_2$.
Given a corrupted image $f\colon\grid\to\MM$, variational models generate a restored image $u\colon\grid\to\MM$
as a minimizer of a functional of the form
\begin{equation} 
\mathcal{J}(u) \coloneqq \mathcal{D}(u;f) +\alpha\mathcal{R}(u),
\end{equation}
where $\mathcal{D}(\cdot;f)$ denotes the data-fitting term, $\alpha>0$ the regularization parameter, and  $\mathcal{R}$ 
the regularization term or prior. 

For real-valued images typically a squared Euclidean distance is chosen as data-fitting term.
For manifold-valued data, we can just use the squared distance on the manifold 
\begin{equation} \label{data_mani}
 \mathcal{D}(u;f) \coloneqq \frac12 \sum_{i\in\grid} \dist(u_i,f_i)^2.
\end{equation}
Depending on the prior knowledge, several regularization terms were proposed in the Euclidean setting.
The (discretized) total variation ($\TV$) introduced by Rudin, Osher, and Fatemi~\cite{ROF92} 
is demonstrably a powerful, edge-preserving, non-smooth, and convex regularizer.
It sums up the norms of the gradients at the image points.
A natural way to define the discrete gradient
\(\nabla  \coloneqq (\nabla_x ,\nabla_y )^\tT\)
on manifold-valued images
is as a vector field in the corresponding tangent spaces
$(T_{u_{i}}\MM)^2$ with the difference operators 
\begin{equation} \label{ictgv:eq:gradUManifoldLog}
	\nabla_x u_i
	\coloneqq
	\begin{cases} \text{if} \;
	\log_{u_{i}}u_{i+(1,0)}\quad& i+(1,0)\in\grid,\\	
	0\quad&\text{otherwise},
	\end{cases}
\end{equation}
and similarly for $\nabla_y$.
Now, the $\TV$ regularizer on manifold-valued images becomes 
\begin{align}\label{eq:tvfunc}
\TV(u) 
&\coloneqq  \sum_{i\in\grid} \biggl( \| \nabla_x u_i\|_{u_i}^p + \nabla_y u_i\|_{u_i}^p  \biggr)^{\frac{1}{p}}  \\
&=\sum_{i\in\grid}\biggl(\sum_{j\in\neighbor(i)}\dist(u_i,u_j)^p\biggr)^{\frac{1}{p}}, \quad p \in \{1,2\},
\end{align}
with $\neighbor(i) \coloneqq  \{i+(1,0),i+(0,1) \}\cap\grid$. Here,
$p=1$ is used for the anisotropic and $p=2$ for the isotropic model. 
This setting was considered in~\cite{LSKC13,WDS14}.
If $\MM = \HH$ is an Hadamard manifold, then functional consisting of the data term \eqref{data_mani} 
and the prior \eqref{eq:tvfunc}
has the same properties as its real-valued version, 
i.e., it is strongly convex and coercive  and hence there exists a unique minimizer.
The TV functional \eqref{eq:tvfunc} is not differentiable. To apply minimization algorithms for differentiable functions
we can recast it, using an even, differentiable function $\varphi: \mathbb R_{\ge 0} \rightarrow \mathbb R_{\ge 0}$,   
in the anisotropic case as
\begin{equation}\label{HalfQuad:functional_1}
\TV_\varphi(u) \coloneqq  \sum_{i \in \grid}  \sum_{j \in \neighbor(i) } \varphi \left( \dist (u_i,u_j) \right),
\end{equation}
and in the isotropic one as 
\begin{equation}\label{HalfQuad:functional_2}
\TV_\varphi(u) \coloneqq  \sum_{i \in \grid} \varphi \biggl( \Bigl( \sum_{j \in \neighbor(i) } \dist (u_i,u_j)^2 \Bigr)^{\tfrac12} \biggr).
\end{equation}
Typical functions $\varphi$ are the Huber function
and $\varphi(x) \coloneqq \sqrt{x^2 + \varepsilon^2}$ with a small $\varepsilon$.

The minimizers of the TV regularized functionals prefer piecewise constant functions,
a behavior called staircasing. To avoid such artifacts second order differences were incorporated
into the regularizer. In the manifold-valued setting, we have to find a counterpart of such differences.
A generalization of the anisotropic so-called second order $\TV$ term was given in \cite{BBSW16}. 
Observing that in the Euclidean case the absolute second order difference
of $x,y,z \in \RR^d$ can be rewritten as
\( \lvert x-2y+z\rvert = 2 \lvert\frac{1}{2}(x+z)-y\rvert\),
a counterpart for
$x,y,z \in \MM$  is defined as
\begin{align}\label{tv:eq:seconddist}
\mathrm{d}_2(x,y,z)
\coloneqq
\min_{c \in\mathcal C_{x,z}}
\dist(c,y),
\end{align}
where \(\mathcal C_{x,z}\) denotes
the set of midpoints~\(\geo{x,z}(\frac12)\)
of all geodesics connecting $x$ and $z$.
Note that the geodesic $\geo{x,z}$ is unique
on Hadamard manifolds.
Similarly,  second order mixed differences were defined for \(x,y,z,w\in\mathcal M\) as
   \[
    \mathrm{d}_{1,1}(x,y,z,w)
    \coloneqq
    \min_{c \in\mathcal C_{x,z}, \tilde c\in\mathcal C_{y,w}}
     \operatorname{dist}(c,\tilde c).
   \]
We emphasize that the absolute second order difference \(\mathrm d_2\) is not convex in \(x\) and $z$
on Hadamard manifolds. 
Now, we can introduce the absolute value of the second order difference in $x$-direction as
     \begin{align}
       \mathrm{d}_{xx} u_{i}
        \coloneqq\begin{cases}
          \mathrm{d}_2(u_{{i} + (1,0)},u_{i},u_{{i} - (1,0)})&\mathrm{if} \; i \pm (1,0) \in \grid, \\
          0&\text{otherwise},
        \end{cases}
     \end{align}
     and similarly in $y$-direction.
As absolute value of the mixed differences we use 
\begin{align}
		\mathrm{d}_{xy} u_{i}
        \coloneqq
    \begin{cases}
      	\mathrm d_{1,1}
          \bigl(
            u_{i},                        u_{i+(0,-1)},
            u_{i+(1,0)},
            u_{i+(1,-1)}
          \bigr)
          & \text{if}
          \ i \pm (0,1) \wedge
      	i + (1,0)\in \grid, \\
      	0 & \mathrm{otherwise},
      	\end{cases}
		\end{align}
    and similarly for $\mathrm{d}_{yx}$. 
		Then we define 
    \begin{align} \label{TV2_mani}
      \operatorname{TV}_2 (u)
        \coloneqq
          \sum_{i \in\grid}
            \Bigl(  \mathrm{d}_{xx}u_{i}^p +\mathrm{d}_{yy}u_{i}^p
        +	    \mathrm{d}_{xy} u_{i}^p  + \mathrm{d}_{xy} u_{i}^p
						\Bigr)^\frac{1}{p},\quad p \in \{1,2\},
      \end{align}
 where $p=1$ is used for the anisotropic model and $p=2$ for the isotropic model.
In the regularizer, the $\TV$ and $ \operatorname{TV}_2$ terms can appear separately or in a coupled way.
Actually, their addition 
$${\mathcal R} (u) \coloneqq \beta \TV(u) + (1-\beta) \operatorname{TV}_2 (u), \quad \beta \in (0,1)$$
was considered in \cite{BBSW16,BW15}. 
Alternatively, couplings which generalize the infimal convolution approach \cite{CL97} to the manifold-valued setting 
were proposed in \cite{BFPS17,BFPS18}.
In the Euclidean setting, the infimal convolution is related to the total generalized variation (TGV) approach 
of Bredies et al.~\cite{BKP10}. Recently, Bredies et al.~\cite{BHSW17} came also up with a TGV model for manifold-valued images,
see also \cite{VBK13} for DT-MRI. In the following we present a different TGV approach from \cite{BFPS18}.
For the relation between both model see \cite[Remark 5.1]{BFPS18}.
In the Euclidean setting, the (discrete) TGV regularizer reads as
$$
\min_{\xi \in (T_u\MM^n)^2}  \Big\{ \sum_{i \in \grid} \beta \lVert \nabla u_i - \xi_i \rVert_{2} + (1-\beta) \lVert  \widetilde \nabla \xi_i \rVert_{2} \Big\},
$$
where $\widetilde \nabla$ denotes a certain symmetric backward difference operator.
For a vector field~\(\xi =(\xi_{i})_{i\in\grid},\ \xi_{i}\in (T_{u_{i}}\MM)^2\), 
the distance between $\xi$ and $\nabla u$ in the first summand is given by
\[
R_1(u,\xi) \coloneqq \sum_{i \in \grid} \left( \lVert \nabla_x u_i - \xi_{1,i} \rVert_{u_i} ^p +  \lVert \nabla_y u_i - \xi_{2,i} \rVert_{u_i} ^p \right) ^\frac1p,
\quad p \in \{ 1,2 \}.
\]
To compute the backward differences $\widetilde \nabla \xi$ in the second summand, 
we need to compare tangent vectors from different tangent spaces.
To this end, we apply the parallel transport between the tangent spaces via the pole ladder.
Since the corresponding expression in \eqref{pre:eq:pole} contains only
exponential and logarithmic maps, the differentials required in the minimization procedure 
can be calculated using the chain rule and Lemma~\ref{pre:differentials}.
More precisely, we define backward differences of a vector field $\zeta \in T_{u}\MM^n$ in $x$-direction by
\begin{equation}
\widetilde{\nabla}_x^{\mathrm{P}} \zeta_i \coloneqq \begin{cases}
\zeta_i - P_{u_{i-(1,0)}\to u_i}^\mathrm{P}(\eta_{i-(1,0)})\quad & \mathrm{if~} i\pm(1,0)\in\grid,\\
0\quad&\mathrm{otherwise,}
\end{cases}
\end{equation}
and similarly in $y$-direction. Then we define
$$
R_2(\xi) 
\coloneqq 
\sum_{i \in \grid} 
\left( \lVert \widetilde{\nabla}_x^{\mathrm{P}} \xi_{1,i} \rVert_{u_i}^p 
+ \lVert \widetilde{\nabla}_y^{\mathrm{P}} \xi_{1,i} \rVert_{u_i}^p 
+  \lVert \widetilde{\nabla}_x^{\mathrm{P}} \xi_{2,i} \rVert_{u_i}^p 
+ \lVert \widetilde{\nabla}_x^{\mathrm{P}} \xi_{2,i} \rVert_{u_i}^p
\right)^\frac1p .
$$
Note that due to simplified computations, this definition differs slightly from the symmetric arrangement of the
backward differences in the Euclidean TGV setting.
Now, we can define a  TGV regularizer for manifold-valued images as 
\begin{equation}\label{ictgv:tgvpole}
\TGV(u)
\coloneqq \inf_{\xi \in (T_u\MM^n)^2} 
\left\{
 \beta R_1(u,\xi) + (1-\beta) R_2(\xi) \right\}, \quad \beta \in (0,1).
\end{equation}
%------------------------------------------------------------------------------
%-----------------------------------------------------------------------------
\section{Minimization Algorithms} \label{sec:algs}
%-----------------------------------------------------------------------------
To compute a minimizer of our functionals, Riemannian optimization methods can be applied.
These intrinsic methods 
are often very efficient 
since they exploit the underlying geometric structure of the manifold, 
see e.g.~\cite{AMS08,RW12}.
For smooth functionals, various methods have been proposed, 
reaching from simple gradient descents on manifolds to more sophisticated trust region or (quasi) Newton methods.
We start by recalling the gradient decent algorithm or more precisely the subgradient algorithm 
which can also be applied for the minimization of non-differentiable functions.

%--------------------------------------------------------
\subsection{Subgradient Descent}
%--------------------------------------------------------
The subdifferential of a convex function $\varphi\colon \MM \to (-\infty,+\infty]$
at $x \in \dom f$ is defined by 
\[
\partial \varphi (x)
\coloneqq
\bigl\{
v \in T_x\MM\colon \varphi(y) \ge \varphi(x) 
+ \langle v,\dotgeodesic{x,y}(0) \rangle
\text{ for all } y \in \dom \varphi
\bigr\},
\]
see, e.g.,~\cite{LMWY11} or~\cite{Udr94} for finite functions $\varphi$.
For any $x \in \interior (\dom \varphi)$, the subdifferential is a nonempty convex and compact set in $T_x\MM$.
If the Riemannian gradient $\grad_\MM \varphi(x)$ of $\varphi$ in $x \in \MM$ exists, 
then 
$\partial \varphi (x) = \bigl\{\grad_\MM \varphi(x)\bigr\}$. 
Further, we see from the definition that $x \in \MM$ 
is a global minimizer of $\varphi$ if and only if $0 \in \partial \varphi (x)$.\\
Let $\varphi\colon\MM\to\RR$ be a convex function and $x^{(0)}\in\MM$ the starting point. Given a sequence $(\tau_r)_{r\in\N}$ of nonegative numbers, 
the subgradient algorithm iterates
\begin{align}
\mbox{for} \; &r=0,1,\ldots \; \mbox{until a stopping criterion is reached}\\
&s^{(r+1)} \in\partial \varphi(x^{(r)}),\\
&x^{(r+1)} =\exp_{x^{(r)}}\Bigl(-\tau_r \frac{s^{(r+1)} }{\lVert s^{(r+1)}\rVert_{x^{(r)}}}\Bigr).
\end{align}

We have the following convergence result, see~\cite{FO98}. 

\begin{theorem}[Convergence of subgradient algorithm]
	Let $\MM$ be a Riemannian manifold with non-negative sectional curvature, $\varphi\colon \MM \to \RR\cup\{\infty\}$
	a convex function which has a minimizer, and
	$(\tau_r)_{r\in\N}$  a sequence of positive numbers in $\ell^2 \backslash \ell^1$.
	Then the sequence $\{x^{(r)}\}_{r\in\N}$ generated by subgradient algorithm converges to a minimizer of $\varphi$.
\end{theorem}

For manifolds with curvature bounded from below the subgradient algorithm converges if the iterates stay in bounded sets, 
see \cite{Wang17} or \cite{WLY15}.

%-------------------------------------------------------------------------------
\subsection{Half-Quadratic Minimization}
%-------------------------------------------------------------------------------
Half-quadratic minimization methods
belonging to the group of quasi-Newton methods are efficient minimization algorithms 
for the functionals with differentiable 
anisotropic or isotropic regularizers $\TV_\varphi$.
These methods, which cover iteratively re-weighted least squares methods, 
were recently generalized to manifold-valued images~\cite{BCHPS16,GS14}.
There exist additive and multiplicative versions of the method, see \cite{GR92,GY95,NN05}.
Here, we focus on the multiplicative one for the isotropic $\TV_\varphi$ regularizer, i.e.,
we want to minimize
$$
{\mathcal J}_\varphi(u) \coloneqq \frac12 \sum_{i\in\grid} \dist(u_i,f_i)^2 
+ \alpha \sum_{i \in \grid} \varphi \biggl( \Bigl( \sum_{j \in \neighbor(i)^+ } \dist (u_i,u_j)^2 \Bigr)^{\tfrac12} \biggr).
$$
We consider the method based on the so-called $c$-transform. 
Given a function $c\colon \RR \times \RR \rightarrow \RR$,
the $c$-transform of a function $\varphi\colon \RR \rightarrow \RR$
is defined by
\[
\varphi^c (s) \coloneqq \inf_{t \in \RR} \bigl\{ c(t,s) - \varphi(t) \bigr\}.
\]
We see immediately that
$
\varphi(t) + \varphi^c (s) \le c(t,s)
$.
For $c(t,s) \coloneqq -st$, the function $\varphi^c = -(-\varphi)^*$ is just the Fenchel transform of $\varphi$.
We need the following proposition, see \cite{BCHPS16}.

\begin{proposition} \label{HalfQuad:prop_1}
Let $\varphi\colon \RR \rightarrow \RR_{\ge 0}$ be an even, differentiable function
and  $c(t,s)\coloneqq  t^2 s$.
\begin{enumerate}
\item[i)] If the function
$\Phi (t) \coloneqq  -\varphi(\sqrt{t})$ for $t \ge 0$ and $\Phi (t)\coloneqq +\infty$ for $t < 0$
is convex, then
$
\varphi = \varphi^{cc},
$
i.e., for $\psi(s) \coloneqq \varphi^c(s)$ it holds
\begin{align} 
  \varphi(t) = \inf_{s \in \mathbb R} \bigl\{c(t,s) - \psi(s) \bigr\} \label{HalfQuad:dual_2}.
\end{align}
\item[ii)] If in addition 
$
\lim_{t \rightarrow \infty} \frac{\varphi(t)}{t^2} \rightarrow 0,
$
and $\varphi'(t) \ge 0$ for $t \ge 0$ and $\varphi''(0+) \coloneqq \lim_{t \rightarrow 0+} \frac{\varphi'(t)}{t}$ exists,
then the infimum in \eqref{HalfQuad:dual_2}
is attained for the tuple $(t,s) = \bigl(t,s(t)\bigr)$ with
\begin{align} \label{HalfQuad:mini_1}
s(t) &\coloneqq 
\begin{cases}
  \frac{\varphi'(t)}{2t} &\text{for } t > 0,\\[1ex]
  \frac{\varphi''(0+)}{2}&\text{for } t = 0.
\end{cases}
\end{align}
These pairs fulfill $\varphi(t) + \psi(s) = c(t,s)$.
The choice is unique except for $t=0$, where any $s$ larger than $\frac{\varphi''(0+)}{2}$
is also a solution.
\item[iii)] If  $\varphi'(t) > 0$ for $t > 0$ and
$\varphi''(0+) > 0$, then $s(t) \in (0,\frac{\varphi''(0+)}{2}]$ for all $t>0$.
\end{enumerate}
\end{proposition}

Functions $\varphi$ fulfilling the assumptions of Proposition \ref{HalfQuad:prop_1} are listed in Table \ref{HalfQuad:possible_phi}.

{\small
\begin{table}
	\setlength{\tabcolsep}{1em}
	\centering
	\begin{tabular}{lll}\toprule
		&$\varphi(t)$ & $s(t)$ 
		\\\midrule
		$\varphi_1(t)$ & $\sqrt{t^2+\varepsilon^2}$ & $\frac{1}{2\sqrt{t^2+\varepsilon^2}}$ \\[1ex]
		$\varphi_2(t)$ & $\begin{cases}\frac{1}{2}t^2\quad &t<\varepsilon,\\
		\varepsilon\lvert t\rvert-\frac{1}{2}\varepsilon^2\quad& t\ge\varepsilon \end{cases}$ & 
		$\begin{cases}\frac{1}{2}\quad &t<\varepsilon,\\
		\frac{\varepsilon}{2\lvert t\rvert}\quad& t\ge\varepsilon 
		\end{cases}$\\[2.5ex]
		$\varphi_3(t)$ & $1-\exp(-\varepsilon^2 t^2)$ & $\varepsilon^2 \exp(-\varepsilon^2 t^2)$\\
		\bottomrule
	\end{tabular}
	\caption{Functions $\varphi$ fulfilling the assumptions of Proposition \ref{HalfQuad:prop_1}.}\label{HalfQuad:possible_phi}
\end{table}
}

Then, replacing $\varphi$ in  $\TV_\varphi$ in the isotropic setting in \eqref{HalfQuad:functional_2} by the expression in \eqref{HalfQuad:dual_2},
we can minimize instead of ${\mathcal J}_\varphi (u)$ the functional
\begin{align} \label{HalfQuad:efunctional_2}
  J(u,v) &\coloneqq 
 \frac12 \sum_{i \in \grid} \dist^2(u_i,f_i) 
 + \alpha  \sum_{i \in \grid}
  \left( c  (\mathrm{d}_i , v_{i} ) - \psi(v_{i}) \right),
\end{align}
where
$\mathrm{d}_i \coloneqq \bigl(\sum_{j \in \neighbor(i)^+ } \dist^2(u_i,u_j) \bigr)^{\frac12}$. 
We  apply alternating minimization over  $v = (v_i)_{i \in \grid} \in \RR^n$ and $u \in \MM^n$ and obtain together with
\eqref{HalfQuad:mini_1} the following iterations:
\begin{align} 
\mbox{for} \; &r=0,1,\ldots \; \mbox{until a stopping criterion is reached}\\
  &v^{(r+1)} =  \argmin_v  J(u^{(r)},v)  =  s \bigl(\mathrm{d}^{(r)} \bigr),\\
  &u^{(r+1)} \in \argmin_u   J(u,v^{(r+1)}). \label{alg:hq}
\end{align}
The minimization over $u$ means  to find a minimizer of
\begin{equation}\label{HalfQuad:efunctional_2_u}
  J (u,v^{(r)}) \coloneqq 
  \frac12 \sum_{i \in \grid } \dist^2(u_i,f_i) + \alpha  \sum_{i \in \grid}
 \Bigl( \sum_{j \in \neighbor(i)^+ } \dist^2(u_i,u_j) \Bigr) v_{i}^{(r)}.
\end{equation}
Here, we can apply, e.g., a gradient descent or a Riemann-Newton method, see~\cite{AMS08}.
Concerning the convergence of the algorithm we have the following theorem, see \cite{BCHPS16}.
%
%-----------------------------------------------------------------------
\begin{theorem}[Convergence of half-quadratic minimization] \label{HalfQuad:convergence}
Let $\HH$ be an Hadamard manifold and 
let $\varphi\colon \mathbb R \rightarrow \mathbb R_{\ge 0}$ fulfill the assumptions of Proposition \ref{HalfQuad:prop_1}.
Then, the sequence $\{u^{(r)} \}_{r \in \mathbb N}$ generated by \eqref{alg:hq}
converges to the uniquely determined minimizer of ${\mathcal J}_\varphi$.
\end{theorem}
%---------------------------------------------------------------------------------------------
\subsection{Proximal Point and Douglas-Rachford Algorithm}\label{sec:prox}
%---------------------------------------------------------------------------------------------
%
In the Euclidean setting,  tools from convex analysis, in particular powerful algorithms 
based on duality theory, were successfully applied  to minimize the proposed functionals.
A prominent example is the alternating
directions method of multipliers (ADMM), which is equivalent to the Douglas-Rachford
algorithm. A central ingredient of these algorithms are proximal mappings which can be efficiently
computed for special regularization terms appearing in Euclidean image processing tasks,
see \cite{BPCPE11,BSS2016}.
Recently, several attempts have been made to translate these concepts to manifolds
and it turns out that on Hadamard manifolds a certain theory of convex functions
can be established. 
For example, the (inexact) cyclic proximal point algorithm can be introduced
on these manifolds \cite{Bac14}, and this method was also used to minimize the functional with first and
second order $\TV$ regularizers in \cite{BBSW16,BLSW14,BHSW17,WDS14}. 
Since the classical Douglas-Rachford algorithm relies on point reflections, it
was natural to extend this algorithm to symmetric Hadamard manifolds \cite{BPS16}.

%------------------------------------------------------------------------------
\subsubsection{Proximal Mapping}
%------------------------------------------------------------------------------
For $\lambda>0$ and a proper, convex, lsc function
$\varphi\colon \mathbb R^m \rightarrow (-\infty,+\infty]$, 
the proximal mapping $\prox_{\lambda \varphi}: \mathbb R^d \rightarrow \mathbb R^d$ defined by
\begin{equation} \label{prox_hilbert}
	\prox_{\lambda \varphi}(x)
	\coloneqq \argmin_{y \in \mathbb R^d} \Big\{ \frac{1}{2} \lVert x-y\rVert_2^2 + \lambda\varphi (y) \Big\}
\end{equation}
is uniquely determined. 
The counterpart  on manifolds reads 
for $\varphi:\MM^m \rightarrow [-\infty,+\infty]$ as
\begin{equation} \label{eq:prox_mani}
\prox_{\lambda \varphi }(x)
\coloneqq 
\argmin_{y \in {\MM^d}} \Big\{  \frac{1}{2}  \dist (x,y) ^2
	+ \lambda\varphi(y)\Big\}.
\end{equation} 
Indeed, for proper, convex, lower semi-continuous  functions on Hadamard manifolds $\MM = \HH$,
the above minimizer is uniquely determined~\cite{Jost97}.
Moreover, the proximal operator is nonexpansive.
For the distance functions appearing in the sums of our models
$\varphi_0^p (x) \coloneqq  \dist(x,y)^p$, $p \in \{1,2\}$ and
$\varphi_1^p (x,y) \coloneqq  \dist(x,y)^p$, $p \in \{1,2\}$,
the proximal mapping can be given analytically, see~\cite{FO02,WDS14}.
For our absolute second order differences
$\varphi_2 (x,y,z) \coloneqq  {\rm d}_2\big(x,y,z\big)$ 
the proximal mapping can be computed numerically on certain manifolds by the (sub)gradient descent algorithm and
Lemma \ref{pre:differentials} as outlined in \cite{BBSW16}.
An analytical expression for the proximal mapping of $\varphi_2$ on the sphere $\mathbb S^1$ was given in \cite{BLSW14}.

The results are summarized in the following lemmas.

%-------------------------------------
\begin{proposition}[Proximal mapping of distance functions] \label{lem:proxies}
	Let $\HH$ be an Hadamard manifold, $\lambda > 0$ and $y \in \HH$.
	\begin{itemize}
		\item[i)]
		The proximal mappings of $\varphi_0^p(x) \coloneqq \frac{1}{p} \dist(\cdot,y)^p$, $p \in \{1,2\}$,
		are given by 
		\begin{align}
		\prox_{\lambda \varphi_0}(x) = \geodesic{x,y}(\hat t), \qquad 
		\hat t &\coloneqq
		\begin{cases}
		\min\big\{ \tfrac{\lambda}{\dist(x,y)} ,1 \big\}&\mbox{\textup{ if }} p = 1,\\
		\tfrac{\lambda}{1+\lambda} & \mbox{\textup{ if }} p = 2.
		\end{cases}
		\end{align}
		\item[ii)]
		The proximal mappings of
		$\varphi_1^p(x,y) \coloneqq \dist(x,y)^p$, $p \in \{1,2\}$,
		are given by 
		\begin{align}
		 \prox_{\lambda \varphi_1} (x,y)
		&= \bigl(\geodesic{x,y}(\hat t),\geodesic{x,y}(\hat t)\bigr), 
		\qquad
		\hat t \coloneqq
		\begin{cases}
		\min\big\{ \tfrac{\lambda}{\dist(x,y)} , \tfrac{1}{2} \big\}&\mbox{\textup{ if }} \; p = 1,\\
		\tfrac{\lambda}{1+2\lambda}&\mbox{\textup{ if }}p = 2.
		\end{cases}
		\end{align}	
	\end{itemize}
\end{proposition}

To give the analytical expressions for the proximal mappings of $\varphi_\nu$, $\nu \in \{0,1,2\}$,
on $\mathbb S^1$, we represent
its elements by the angles in $[-\pi,\pi)$. For $a \in \mathbb R$, we
denote by $(a)_{2\pi} \in [-\pi,\pi)$ those number for which there exists $k \in \mathbb Z$ such that
$a + 2\pi k = (a)_{2\pi}$.

\begin{proposition}[Proximal mapping of distance functions on $\mathbb S^1$] 
Let $w_1 \coloneqq (-1,1)^\tT$, 
$w_2\coloneqq (1,-2,1)^\tT$ and 
$s_\nu \coloneqq \sgn(\langle x,w_\nu \rangle)_{2 \pi}$, $\nu \in \{1,2\}$.
Then, for $\nu \in \{1,2\}$, the following holds true:
\begin{itemize}
\item[i)]
If $\lvert(\langle x,w_\nu \rangle)_{2\pi} \rvert < \pi$, then 
$$
\prox_{\lambda \varphi_\nu} (x) =  (x - s_\nu \, m_\nu \,w_\nu )_{2\pi}, \qquad
m_\nu \coloneqq \min \left\{\lambda ,\frac{\lvert (\langle x,w_\nu \rangle)_{2\pi} \rvert} {\lVert w_\nu \rVert_2^2} \right\}.
$$
\item[ii)]
If $\lvert(\langle x,w_\nu \rangle)_{2\pi} \rvert = \pi$,
then the proximal mapping is  two-fold
$$
\prox_{\lambda \varphi_\nu} (x) =  (x \pm s_\nu \, m_\nu \, w_\nu)_{2\pi},
 \qquad m_\nu\coloneqq \min\left\{\lambda ,\frac{\pi}{\lVert w_\nu \rVert_2^2} \right\}.
$$
\item[iii)]
If $\lvert(\langle x,w_\nu \rangle)_{2\pi} \rvert < \pi$, then 
\begin{equation}
\prox_{\lambda \varphi_\nu^2} (x) = \left(x-\lambda \, \frac{ (\langle x,w_\nu \rangle)_{2\pi}}{1+\lambda \|w_\nu\|_2^{2} } w_\nu\right)_{2 \pi}.
\end{equation}
\item[iv)]
If $\lvert(\langle x,w_\nu \rangle)_{2\pi} \rvert = \pi$,
then the proximal mapping is  two-fold
\begin{equation}
\prox_{\lambda \varphi_\nu^2} (x) = \left(f \pm \lambda \, \frac{\pi}{1+\lambda \|w_\nu\|_{2}^{2}} w_\nu \right)_{2 \pi}.
\end{equation}
\item[v)] Finally, 
\begin{equation} \label{min_quad_1}
\prox_{\lambda \varphi_0^2} (x) =  \left( \frac{x+\lambda y}{1+\lambda} + \frac{\lambda}{1+\lambda} \, 2\pi \, v \right)_{2\pi},
\end{equation}
where 
\[
v \coloneqq
\left\{
\begin{array}{ll}
0 & \mbox{\textup{ if }} \; \lvert x-y\rvert \le \pi ,\\
\sgn(x - y) & \mbox{\textup{ if }} \; \lvert x-y\rvert > \pi.
 \end{array}
\right.
\]
\end{itemize}
\end{proposition}
%---------------------------------------------------------------------------------------------
\subsubsection{Cyclic Proximal Point Algorithm}\label{sec:cppa}
%---------------------------------------------------------------------------------------------
Our functionals have the general form
\begin{equation}\label{eq:sumoffunc}
{\mathcal J}(x)\coloneqq \sum_{k=1}^K \varphi_k(x),
\end{equation} 
where appropriate splittings into the $K$ summands must be determined.
Given a starting point $x^{(0)}\in\MM$ the cyclic proximal point algorithm (CPPA) iterates
\begin{align}\label{eq:cppa}
\mbox{for} \; &r=0,1,\ldots \; \mbox{until a stopping criterion is reached}\\
&\mbox{for} \; k=1,\ldots,K\\
&\quad x^{(r+\frac{k}{K})} \coloneqq \prox_{\tau_r {\varphi}_{k}}   (x^{(r + \frac{k-1}{K})}) .
\end{align}
We have the following convergence result from \cite{BBSW16}.

%-------------------------------------------------------------------------------
\begin{theorem}[Convergence of cyclic PPA] \label{thm:cyclic} 
Let ${\mathcal H}$ be a Hadamard manifold and 
$\varphi_k\colon {\mathcal H} \rightarrow \mathbb R$, $k=1,\ldots,K$, convex continuous functions such that
${\mathcal J}$ attains a (global) minimum.
Assume that there exist $p \in {\mathcal H}$ and $C>0$ such that for each $k=1,\dots,K$ and all $x,y\in {\mathcal H}$ we have 
\begin{equation} \label{i:sppa:lips}
\varphi_k(x)-\varphi_k(y) \leq C \dist(x,y) \left(1+\dist(x,p) \right).
\end{equation}
Then the sequence $\{ x^{(r)} \}_{r \in \mathbb N}$ 
with non-negative $\{\tau_r\}_{r \in \mathbb N} \in \ell^2 \backslash \ell_1$  
converges for every starting point $x^{(0)}$ 
to a minimizer of ${\mathcal J}$.
\end{theorem}

The result can be generalized for the inexact cyclic PPA which iteratively generates the points
$x^{(r+ \tfrac{k}{K})}$, $k=1,\ldots,K$, $r\in \mathbb N_0$,
fulfilling
\begin{equation} \label{inexact_PPA}
\dist\big(x^{(r+ \tfrac{k}{K})},\prox_{\tau_r {\varphi}_k} (x^{( r + \tfrac{k-1}{K} )}) \big) < \frac{\varepsilon_r}{K},
\end{equation}
where $\{ \varepsilon_r \}_{r\in \mathbb N_0}$ is a given sequence of positive reals with
$\sum_{r=1}^\infty \varepsilon_r < \infty$, see \cite{BBSW16}.

%------------------------------------------------------------------------------
\subsubsection{Douglas-Rachford Algorithm for Symmetric Hadamard Spaces}
\label{subsec:DR}
%------------------------------------------------------------------------------
The Douglas-Rachford (DR) algorithm relies on reflections.
In the Euclidean setting, the reflection of a proper, convex, lsc function 
$\varphi\colon \mathbb R^d \rightarrow (-\infty,+\infty]$ is defined as  
\begin{equation}\label{reflection_R_1}
R_\varphi (x) = 2\prox_{ \varphi}(x) - x.
\end{equation}
It is a nonexpansive operator on $\mathbb R^d$ with respect to the Euclidean norm. 
Given two proper, convex, lsc functions
$\varphi,\psi\colon\mathbb R^d \rightarrow (-\infty,+\infty]$, the DR algorithm 
aims to solve
\begin{equation} \label{DR:drs_split_1}
 \argmin_{x \in \mathbb R^d} \bigl\{ \varphi(x) + \psi(x) \bigr\}
\end{equation}
by iterating a starting point $t^{(0)}$ as follows:
%---------------------------------------------------
\begin{align}
\mbox{for} \; &r=0,1,\ldots \; \mbox{until a stopping criterion is reached}\\
 &t^{(r+1)} \coloneqq \bigl(
			(1-\tau_r) \, I
			+ \tau_r R_{\eta \varphi} R_{\eta \psi}
			\bigr) \bigl(t^{(r)}\bigr).\\
	&x^{(r+1)} \coloneqq  \prox_{\eta\psi}(t^{(r+1)}).
	\end{align}
%---------------------------------------------------
%
Note that $x^{(r)}$ must be only computed in the final step of the algorithm.
It is known that the DR algorithm
converges if
$\mathrm{ri}(\dom \varphi) \cap \mathrm{ri}(\dom \psi) \not = \emptyset$, a minimizer exists, $\eta >0$ and
$\sum_{r \in \mathbb N} \tau_r(1-\tau_r) = + \infty$.
The DR algorithm can be considered a special case of the
Krasnoselski--Mann iteration
\begin{equation} \label{DR:it:km_real}
t^{(r+1)} = \bigl( (1-\lambda_r)\, I + \lambda_r T\bigr)(t^{(r)}),
\end{equation}
with $T\coloneqq R_{\eta \varphi} R_{\eta \psi}$.
It is well-known that the sequence of iterates to a fixed point of $T$ if $T$ is nonexpansive. 
This is clearly the case for our setting since 
reflections $R_\varphi$ are nonexpansive and the 
concatenation of nonexpansive functions is nonexpansive again.

For two points \(x,a\in \HH\) on a symmetric Hadamard manifold, 
the geodesic reflection \eqref{def:refl} can be written
as \(R_p(x) = \exp_p(-\log_px)\).
Then, the geodesic reflection of a proper, convex, lsc function
$\varphi\colon\HH^n \rightarrow (-\infty,+\infty]$ is the mapping 
\begin{equation} \label{DR:reflection_hada}
{R}_\varphi (x) = \exp_{\prox_{ \varphi }(x)}\bigl(-\log_{\prox_{ \varphi }(x)}(x)\bigr).
\end{equation}
Now in order to minimize
\begin{equation} \label{DR:drs_split}
 \argmin_{x \in \HH^d} \bigl\{ \varphi(x) + \psi(x) \bigr\}
\end{equation}
the DR algorithm can be generalized as follows:
\begin{align}
\mbox{for} \; &r=0,1,\ldots \; \mbox{until a stopping criterion is reached}\\
		& t^{(r+1)} \coloneqq  \geodesic{t^{(r)},s^{(r)}} (\tau_r), 
		\quad s^{(r)} \coloneqq  {R}_{\eta \varphi} {R}_{\eta \psi} \bigl( t^{(r)} \bigr),\\
		& x^{(r+1)} \coloneqq  \prox_{\eta\psi}(t^{(r+1)}).
	\end{align}
Again, the algorithm can be seen as special case of the Krasnoselski--Mann iteration
\begin{equation}\label{eq:krasnoleski}
 t^{(r+1)}\coloneqq \gamma\bigl(t^{(r)},T(t^{(r)});\tau_r\bigr). 
\end{equation}
It was proved in~\cite{Ka13}, see also~\cite[Theorem 6.2.1]{Bac14b} that such iteration converges
to a fixed point of $T$,
if $T$ is nonexpansive, has a nonempty fixed point set  and $\sum_{r \in \mathbb N} \tau_r(1-\tau_r) = + \infty$.
Unfortunately, reflections at proper, convex, lsc functions on Hadamard manifolds are in general \emph{not nonexpansive}. 
However, for the distance functions involved in our functionals  nonexpansivness is guaranteed
by the following theorem.

%-------------------------------------------------------
\begin{theorem}[Reflections at distance functions] \label{reflex:d(xa)}
For an arbitrary fixed $a \in {\mathcal H}$ and $\varphi(x) \coloneqq \dist^p(a,x)$,
$p \in\{1,2\}$,
the geodesic reflection ${R}_{\eta \varphi}$, $\eta >0$ is nonexpansive.
For $\varphi(x,y)\coloneqq \dist^p (x,y)$, $p \in \{1,2\}$, the geodesic reflection $R_{\eta \varphi}$, $\eta>0$, is nonexpansive.
\end{theorem}
%-------------------------------------------------------

In general we have more than two summands, i.e., we are interested in the minimization of \eqref{eq:sumoffunc}
where $\varphi_k\colon \HH^n \rightarrow (-\infty,+\infty]$, $k=1,\ldots,K,$ are proper, convex, lsc functions.
Here the trick is to rewrite the functional as the sum of two special components
\begin{equation}\label{DR:eq:prod_problem_had}
	\argmin_{x\in\HH^{nK}} \bigl\{\Phi(x)+\iota_{\textsf{D}}(x) \bigr\},
\end{equation}
where 
$\Phi (x)\coloneqq\sum_{k=1}^K \varphi_k(x_k)$, 
$x\coloneqq (x_k)_{k=1}^K,$ 
and
\[
{\textsf{D}}\coloneqq\{x\in\HH^{nK} \colon x_1=\dots=x_n\in\HH^n\}.
\]
Obviously, ${\textsf{D}}$ is a nonempty, closed convex set so that its
indicator function is proper, convex and lsc, see~\cite[p. 37]{Bac14b}.
Now, the DR algorithm can be formulated as
%-------------------------------------
\begin{align}
\mbox{for} \; &r=0,1,\ldots \; \mbox{until a stopping criterion is reached}\\
 &t^{(r+1)} \coloneqq \geodesic{\vec{t}^{(r)},\vec{s}^{(r)}} (\lambda_r), \quad 
  s^{(r)} \coloneqq {R}_{\eta \Phi} {R}_{\iota_{\textsf{D}}}\bigl( \vec{t}^{(r)} \bigr),\\
 &x^{(r+1)} \coloneqq \Pi_{\textsf{D}}(t^{(r+1)}).
\end{align}
%-------------------------------------
Note that the second step is indeed only necessary in the final iteration.
Concerning the last step note that
\begin{equation} \label{DR:prox_iota_D_had}
 \Pi_{\textsf{D}}(x)
  = \Big(
	\argmin_{x \in \HH^n} \sum_{k=1}^K \dist(x_k,x)^2, \ldots, \argmin_{x \in \HH^n} \sum_{k=1}^K \dist(x_k,x)^2
	\Big)
	\in\HH^{nK}	.
 \end{equation}
The minimizer of the sum is the so-called Karcher mean,
which can be efficiently computed on Hadamard manifolds using the gradient descent
algorithm or the cyclic proximal point algorithm.

Concerning the convergence of the parallel DR algorithm, for our setting, ${R}_{\eta \Phi}$ is nonexpansive since it contains only geodesic reflections of the distance functions in Theorem \ref{reflex:d(xa)}.
Unfortunately, in symmetric Hadamard manifolds,
geodesic reflections corresponding to orthogonal projections onto convex sets are in
general not nonexpansive.
This is also true for our special set $\textsf{D}$.
The situation changes if we consider manifolds with constant
curvature $\kappa$. Here those reflections are nonexpansive, see~\cite{FL13,Persch2018}.
So, in summary, although the parallel DR algorithm showed  very good numerical performance on general Hadamard manifolds in \cite{BPS16},
theoretical convergence results remain up to now limited to manifolds with constant non-positive curvature.

%-------------------------------------------------------------------------------------
%-----------------------------------------------------------------------------
\section{Numerical Examples} \label{sec:numerics}
\begin{figure}
\begin{center}
	\begin{tabular}{ccc}
		&\includegraphics[width=0.25\textwidth]{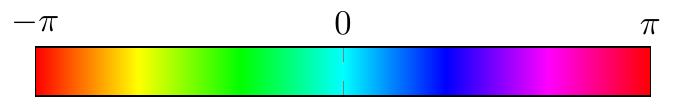}&
		\\[1ex]
		\includegraphics[width=0.25\textwidth]{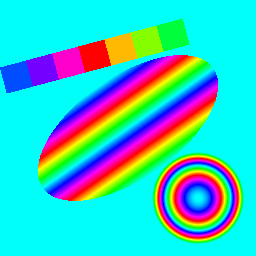}&
		\includegraphics[width=0.25\textwidth]{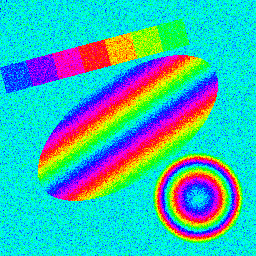}&
		\includegraphics[width=0.25\textwidth]{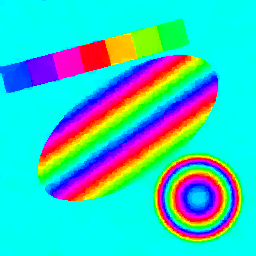}
		\\
		 Original&   Noisy image ($88.5 \times 10^{-3}$)& TV  ($7.2 \times 10^{-3}$)
		\\
		\includegraphics[width=0.25\textwidth]{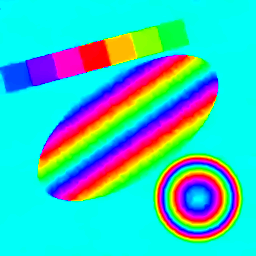}&
		\includegraphics[width = 0.25\textwidth]{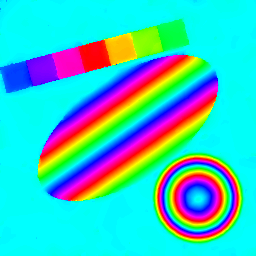}&
		\includegraphics[width=0.25\textwidth]{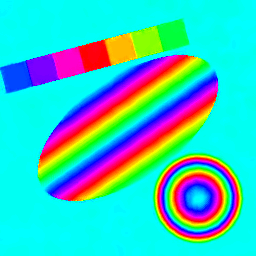}
		\\
		 TV-TV2  ($5.2 \times 10^{-3}$)& TGV  ($2.6 \times 10^{-3}$)&NL-MMSE ($2.5 \times 10^{-3}$)		
	\end{tabular}
	\end{center}
	\caption{Comparison of different variational models for an image with values on $\mathbb S^1$.} \label{fig:circle} 
\end{figure}

\begin{figure}
	\centering
	\begin{tabular}{ccc}
		\includegraphics[width = 0.3\textwidth]{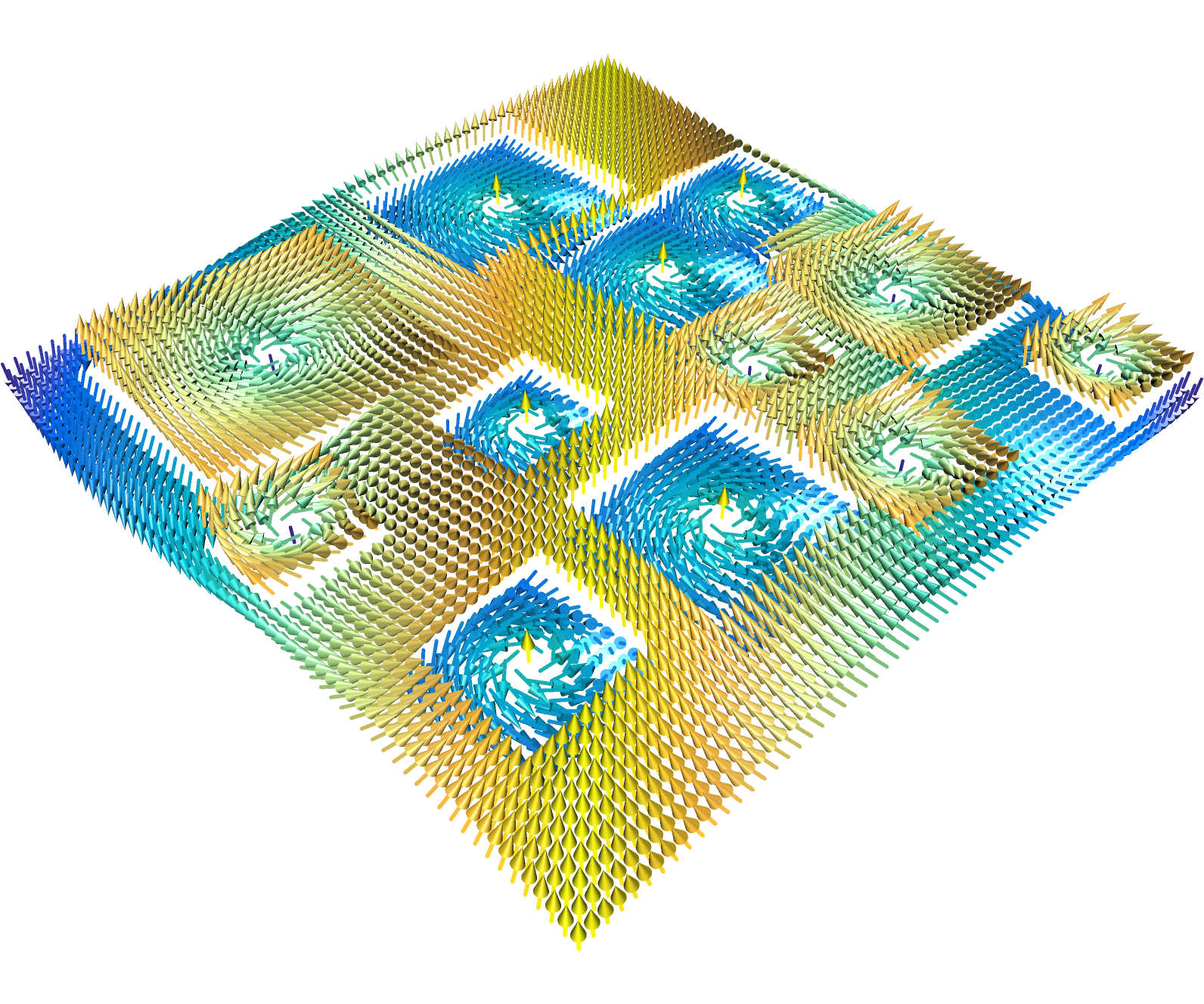}&
		\includegraphics[width = 0.3\textwidth]{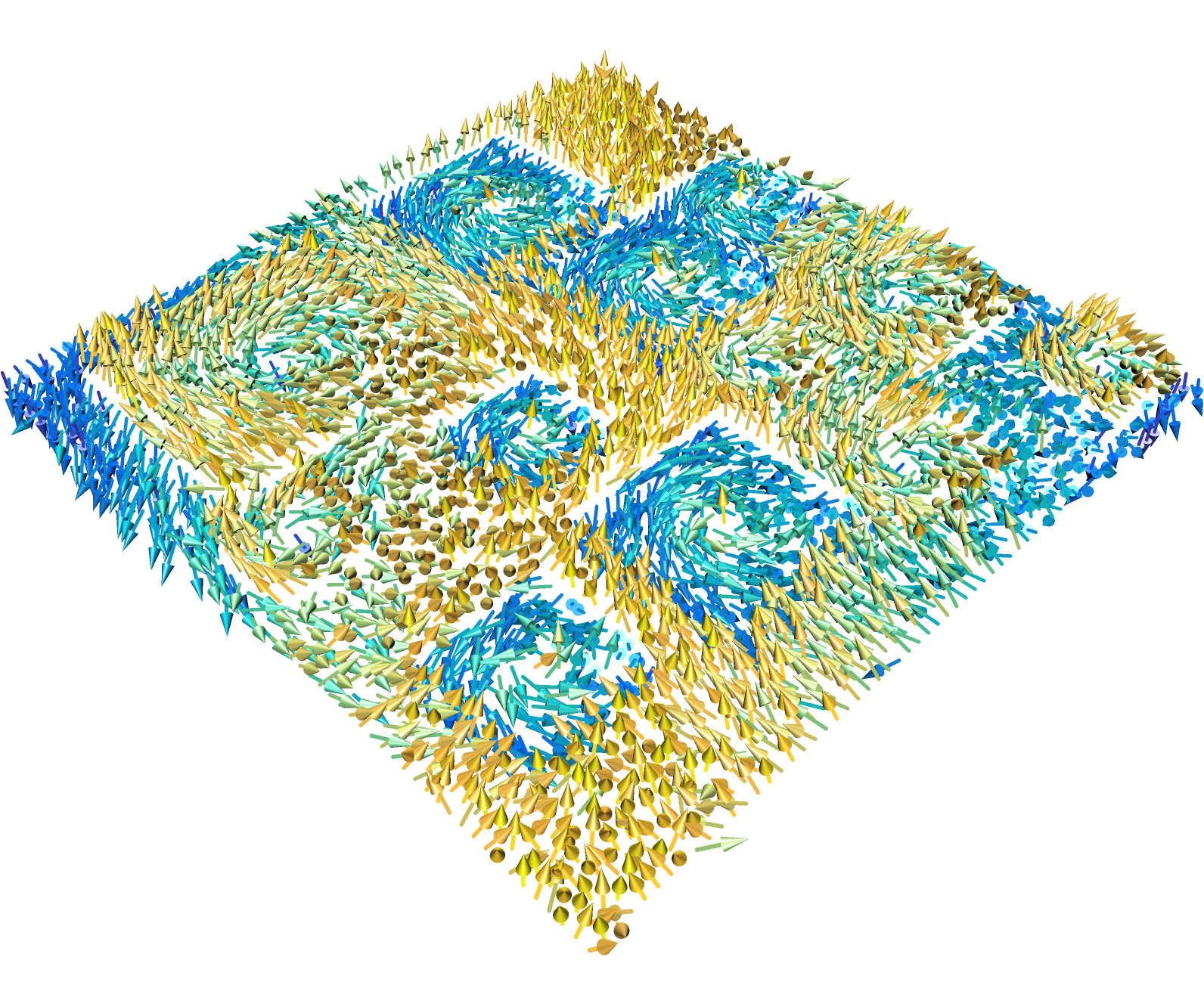}&
		\includegraphics[width = 0.3\textwidth]{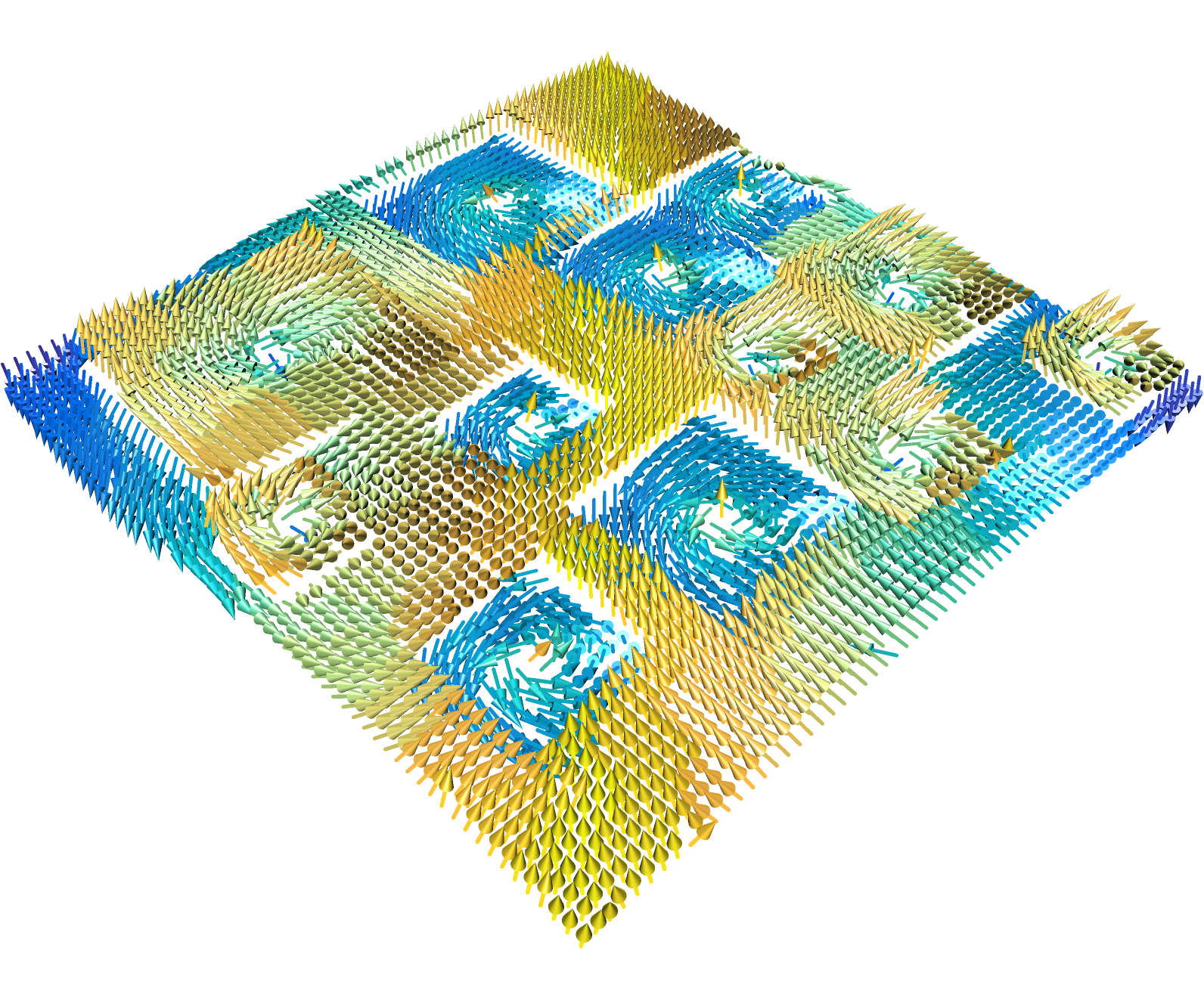}\\
		 Original image&Noisy image ($0.1767$)&TV  $(0.0352)$\\
		 \includegraphics[width = 0.3\textwidth]{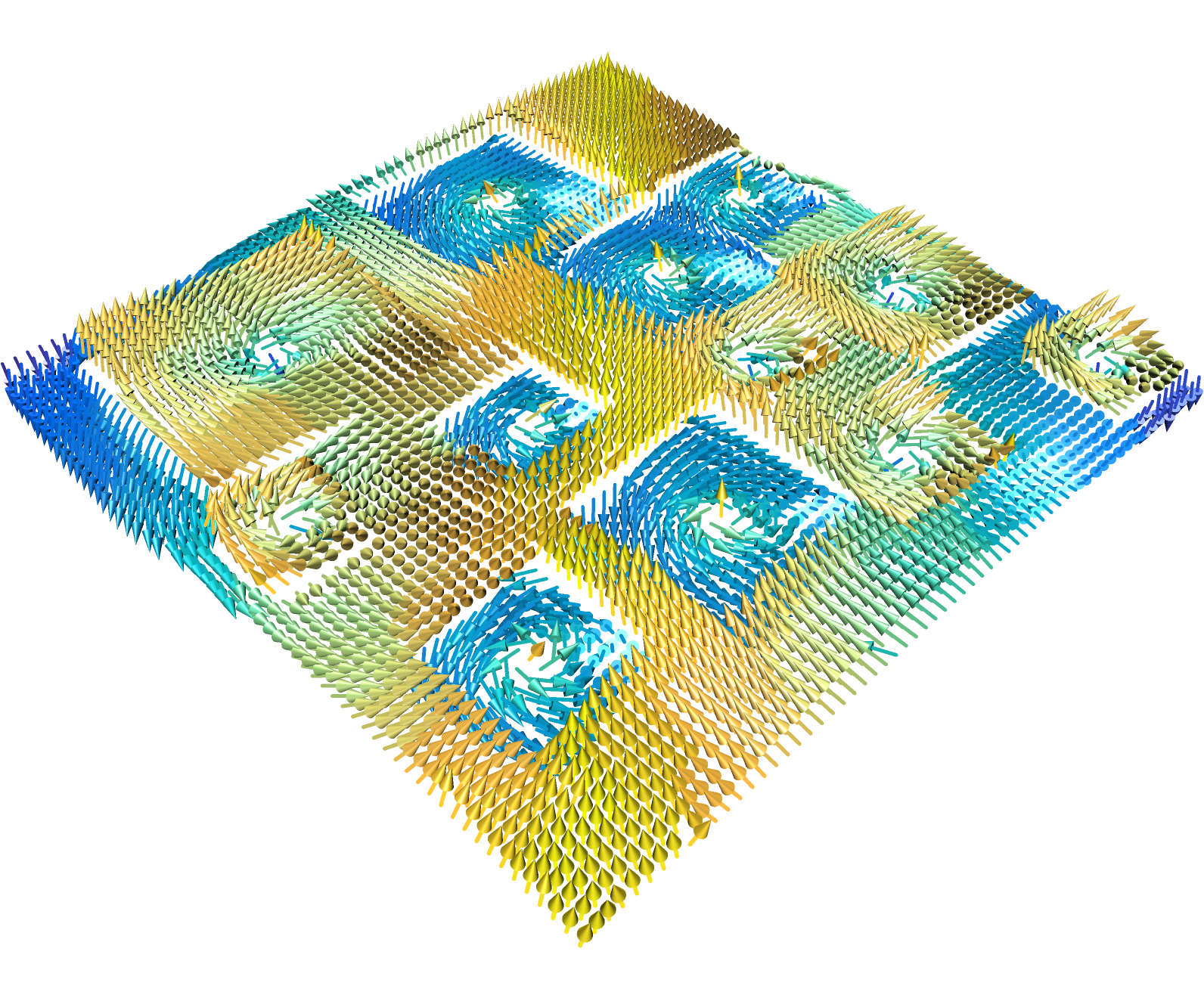}&
		\includegraphics[width = 0.3\textwidth]{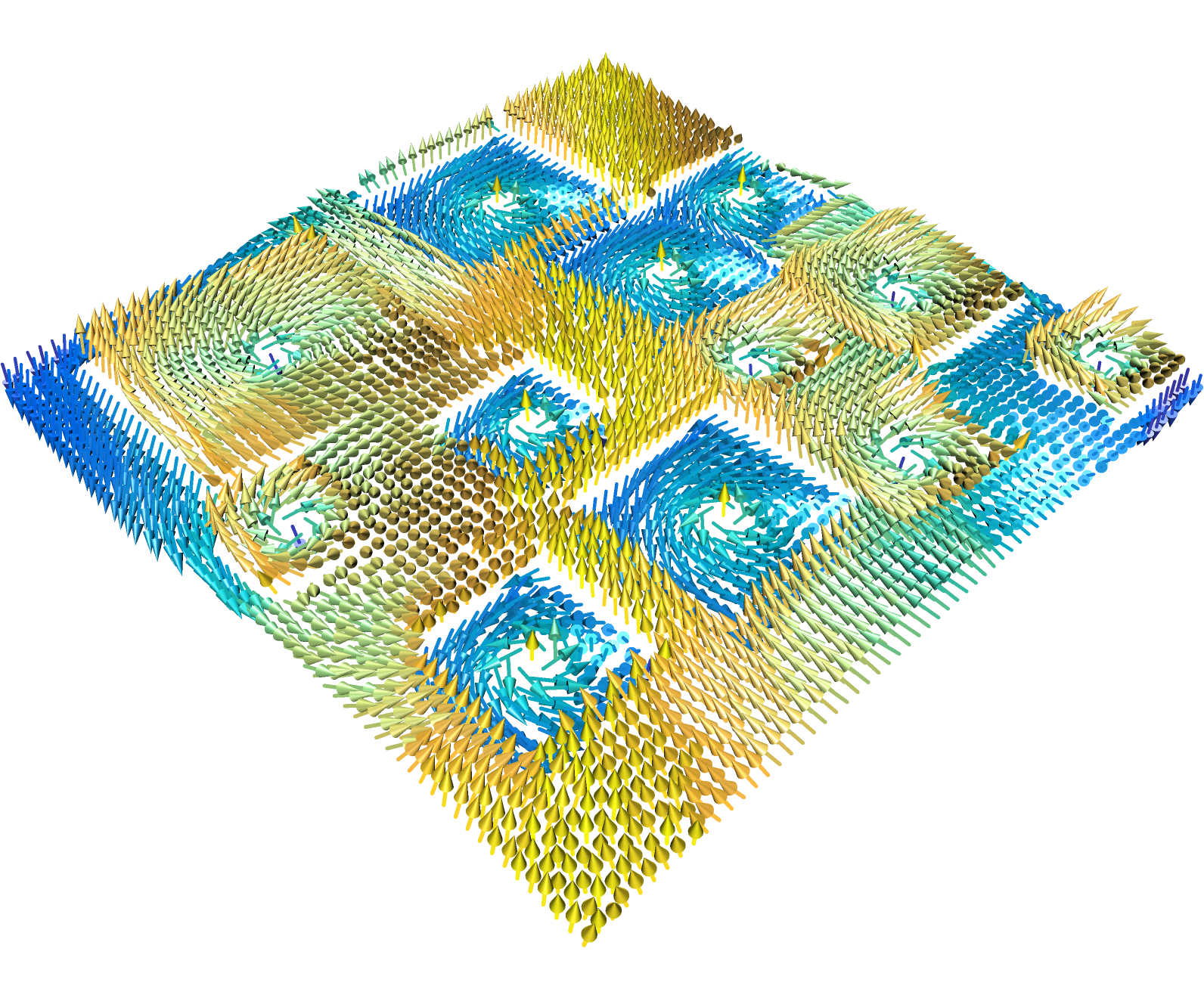}&
		\includegraphics[width = 0.3\textwidth]{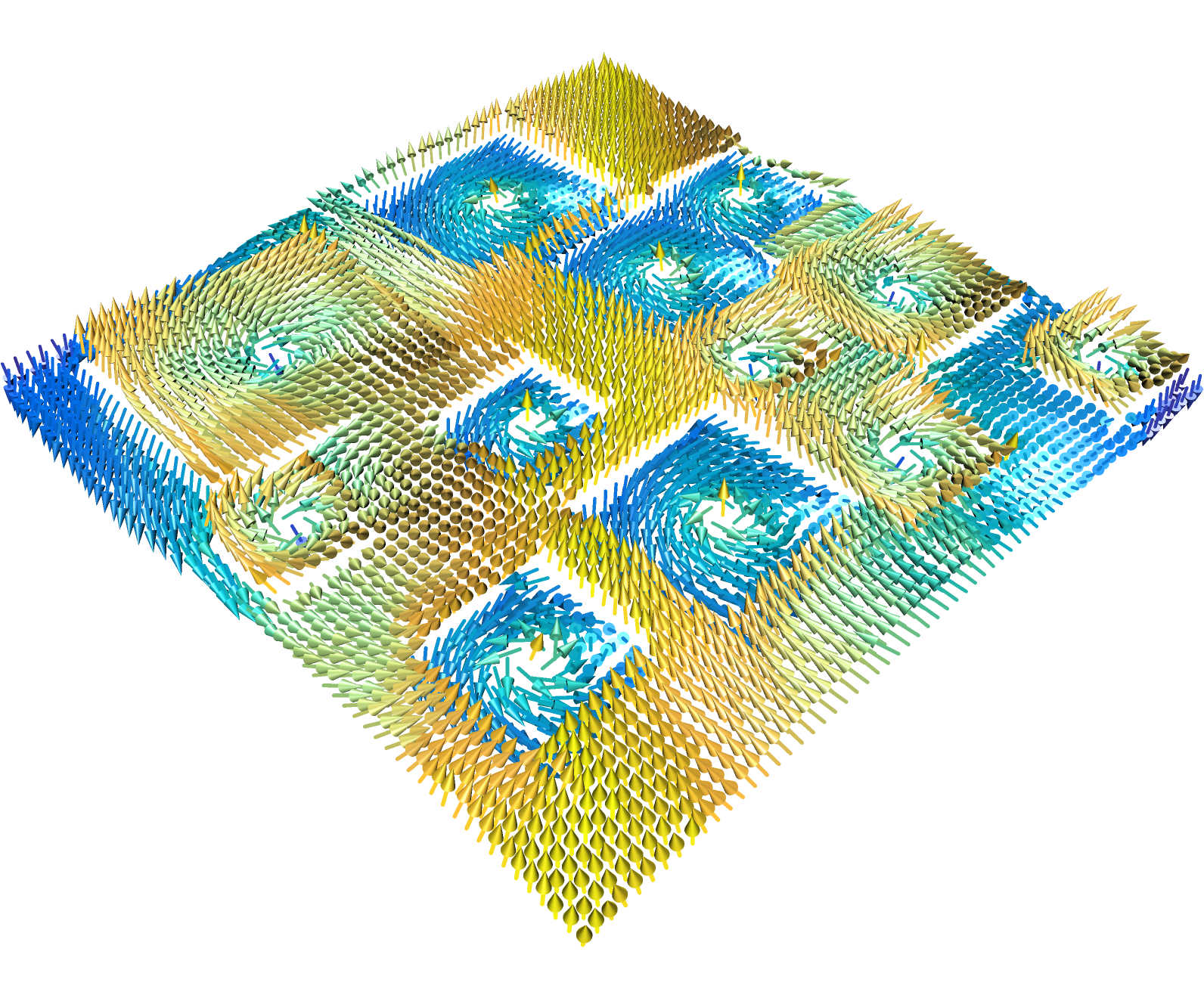}\\
		TV-$\operatorname{TV}_2$  ($0.0338$) &	NL-means ($0.0326$) &	NL-MMSE ($0.0258$)
	\end{tabular}	
	\caption{Comparison of different variational models for an image with values on $\mathbb S^2$.} \label{fig:sphere} 
\end{figure}

%-----------------------------------------------------------------------------
In this section, we give some illustrative numerical examples.
The experiments are carried out   using \textsc{Matlab} 2017a
and the MVIRT toolbox~\cite{Bergmann2017}\footnote{Open source, available at \href{https://ronnybergmann.net/mvirt}{ronnybergmann.net/mvirt/}}.
As a quality measure we use the
mean squared error (MSE) defined by
$$
\epsilon\coloneqq
\tfrac{1}{\lvert\grid\rvert}
\sum_{{i} \in\grid}\dist^2(u_{i},u_{0,{i}}),
$$
where $u_0$ denotes the original image.
The parameters in the models were obtained via a grid search with respect to the optimal
$\epsilon$ and can be found in detail in the respective papers \cite{BCHPS16,BFPS18}.
In Figures~\ref{fig:circle} and~\ref{fig:sphere} we compare the performance of our variational models with different regularizers,
where TV-TV2 denotes the additive coupling of $\TV$ and $\TV_2$.
For comparison we added the results obtained with the patched-based methods from \cite{BLPS17}, namely with nonlocal means (NL-means)
and nonlocal MMSE (NL-MMSE). In brackets we give the corresponding error value $\epsilon$.

\begin{figure}
	\centering
	\begin{tabular}{ccc}
		\includegraphics[width = 0.3\textwidth]{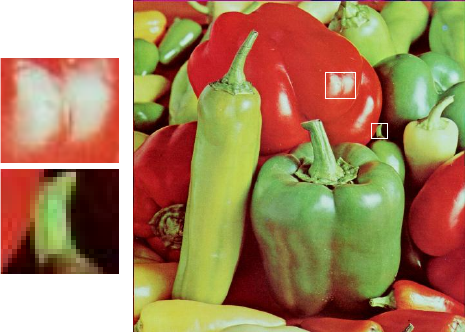}&
		\includegraphics[width = 0.3\textwidth]{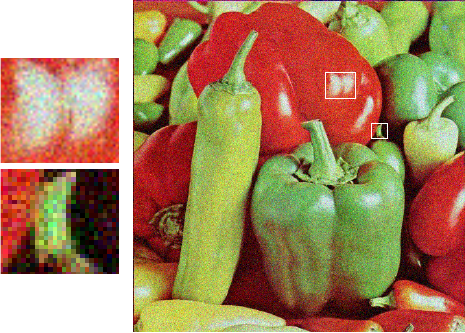}&
		\includegraphics[width = 0.3\textwidth]{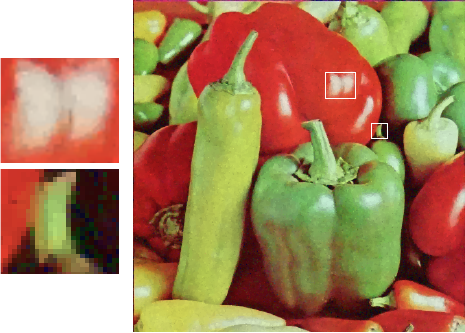}\\
		 Original image & Noisy image (20.31)&(nonsmoothed) TV  (29.32)\\[0.5ex]
		 \includegraphics[width = 0.3\textwidth]{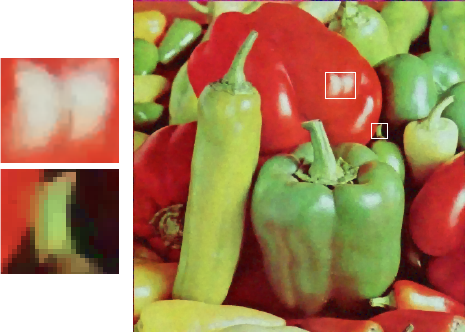}&
		\includegraphics[width = 0.3\textwidth]{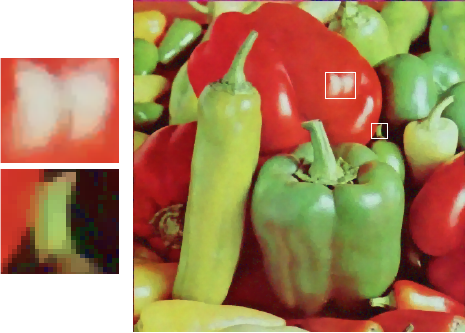}&
		\includegraphics[width = 0.3\textwidth]{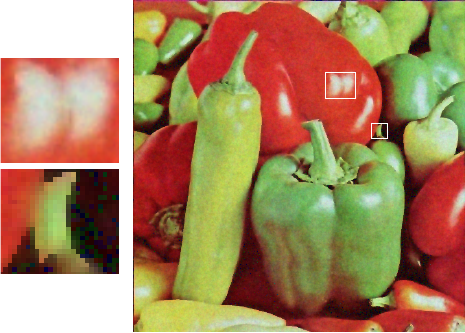}\\
		HQ with~$\varphi_1$  (30.29)
    & HQ with~$\varphi_2$ (29.68)
    &	HQ with~$\varphi_3$ (28.95)
	\end{tabular}	
	\caption{Comparison of  half-quadratic minimization applied to the anisotropic $L_2$-$\TV$
	model with different functions $\varphi_i$, $i=1,2,3$,
	on the chromaticity-brightness color model, 
	i.e.~the manifold $\mathbb R\times\mathbb S^2$.} \label{fig:HQ} 
\end{figure}

\begin{figure}[tbp]
\centering
				\includegraphics[width=.3\textwidth]{Images_Book/CaminoOrig} 
				\hspace{0.5cm}
				\includegraphics[width=.3\textwidth]{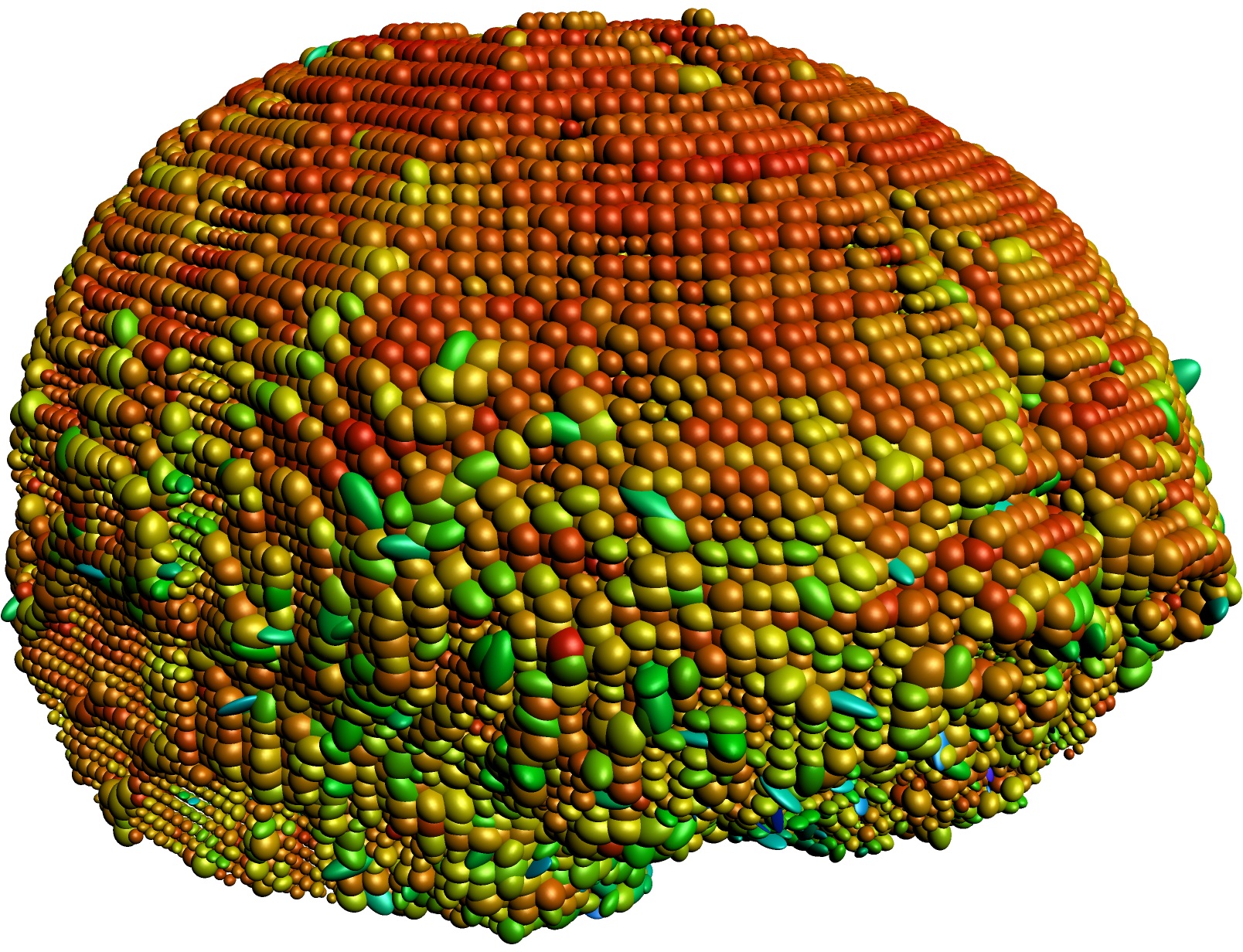}
		\caption{Original ``Camino'' data set (left) and a denoised version by the smoothed $\TV$ regularizer with $\varphi_1$ (right).}\label{fig:DTMRI}
\end{figure}

Figures~\ref{fig:HQ} and~\ref{fig:DTMRI} show denoising results obtain by the half quadratic minimization.
Figures~\ref{fig:HQ} contains results for the different functions $\varphi$ in Table~\ref{HalfQuad:possible_phi}
and for the nonsmooth $\TV$ regularizer. Here the parameters are optimized with respect to the PSNR of the RGB images
and the numbers in brackets refer to the PSNR.
In Figure~\ref{fig:DTMRI} we present a denoising result for the 3D DT-MRI image from Figure~\ref{fig:ex}.

%------------------------------------------------------------------------------------
%-------------------------------------------------------------------------------

%-----------------------------------------------------------------------------
%\section{CONCLUSIONS} \label{sec:conclusions}
%-----------------------------------------------------------------------------

\bibliographystyle{abbrv}
\bibliography{dis}

\end{document}